# Exponential bounds for minimum contrast estimators

**Yuri Golubev**

*Université de Provence*
*39, rue F. Joliot-Curie*
*13453 Marseille, France*
*e-mail:* `golubev.yuri@gmail.com`

**Vladimir Spokoiny**

*Weierstrass-Institute and*
*Humboldt University Berlin,*
*Mohrenstr. 39, 10117 Berlin, Germany*
*e-mail:* `spokoiny@wias-berlin.de`

**Abstract:** The paper focuses on general properties of parametric minimum contrast estimators. The quality of estimation is measured in terms of the rate function related to the contrast, thus allowing to derive exponential risk bounds invariant with respect to the detailed probabilistic structure of the model. This approach works well for small or moderate samples and covers the case of a misspecified parametric model. Another important feature of the presented bounds is that they may be used in the case when the parametric set is unbounded and non-compact. These bounds do not rely on the entropy or covering numbers and can be easily computed. The most important statistical fact resulting from the exponential bonds is a concentration inequality which claims that minimum contrast estimators concentrate with a large probability on the level set of the rate function. In typical situations, every such set is a root-n neighborhood of the parameter of interest. We also show that the obtained bounds can help for bounding the estimation risk, constructing confidence sets for the underlying parameters. Our general results are illustrated for the case of an i.i.d. sample. We also consider several popular examples including least absolute deviation estimation and the problem of estimating the location of a change point. What we obtain in these examples slightly differs from the usual asymptotic results presented in statistical literature. This difference is due to the unboundness of the parameter set and a possible model misspecification.

## 1. Introduction

One of the most fundamental ideas in statistics is to describe an unknown distribution $I\!P$ of the observed data $\boldsymbol{Y} \in I\!R^n$ with the help of a simple parametric family $(I\!P_{\boldsymbol{\theta}}, \boldsymbol{\theta} \in \Theta)$, where $\Theta$ is a subset in a finite dimensional space, say, in $I\!R^p$. In this situation, the statistical model is characterized by the value of the parameter $\boldsymbol{\theta} \in \Theta$ and the statistical inference about $\mathbb{P}$ is reduced to recovering $\boldsymbol{\theta}$. The standard likelihood approach suggests to estimate $\boldsymbol{\theta}$ by maximizing the corresponding likelihood function. The maximum likelihood estimator can be generalized in several ways resulting in the so-called *minimum contrast* and *M-estimators*; see Huber (1967) and Huber (1981). The main idea behind this generalization is to estimate the underlying parameter $\boldsymbol{\theta}$ by minimizing over $\Theta$ a *contrast function* $-L(\boldsymbol{Y}, \boldsymbol{\theta})$:

$$\widetilde{\boldsymbol{\theta}} = \underset{\boldsymbol{\theta} \in \Theta}{\operatorname{argmin}}\{-L(\boldsymbol{Y}, \boldsymbol{\theta})\} = \underset{\boldsymbol{\theta} \in \Theta}{\operatorname{argmax}} L(\boldsymbol{Y}, \boldsymbol{\theta}). \tag{1.1}$$

The negative sign in this notation comes from the main example which we have in mind when $L(\boldsymbol{Y}, \boldsymbol{\theta})$ is the log-likelihood or quasi log-likelihood. A natural condition on the contrast function is that its expectation under the true measure $\mathbb{P}_{\boldsymbol{\theta}_0}$ is minimized at the true parameter $\boldsymbol{\theta}_0$, i.e.

$$\boldsymbol{\theta}_0 = \underset{\boldsymbol{\theta} \in \Theta}{\operatorname{argmax}} \mathbb{E}_{\boldsymbol{\theta}_0} L(\boldsymbol{Y}, \boldsymbol{\theta}). \tag{1.2}$$

If $L(\boldsymbol{Y}, \boldsymbol{\theta})$ is log-likelihood ratio, that is,

$$L(\boldsymbol{Y}, \boldsymbol{\theta}) = \log \frac{d\mathbb{P}_{\boldsymbol{\theta}}}{d\mathbb{P}_{\boldsymbol{\theta}_0}}(\boldsymbol{Y})$$

then the value $-\mathbb{E}_{\boldsymbol{\theta}_0} L(\boldsymbol{\theta}, \boldsymbol{\theta}_0)$ coincides with the Kullback-Leibler divergence $\mathcal{K}(\mathbb{P}_{\boldsymbol{\theta}_0}, \mathbb{P}_{\boldsymbol{\theta}})$ between $\mathbb{P}_{\boldsymbol{\theta}_0}$ and $\mathbb{P}_{\boldsymbol{\theta}}$. It is well known that $\mathcal{K}(\mathbb{P}_{\boldsymbol{\theta}_0}, \mathbb{P}_{\boldsymbol{\theta}})$ is always non-negative and $\mathcal{K}(\mathbb{P}_{\boldsymbol{\theta}_0}, \mathbb{P}_{\boldsymbol{\theta}}) = 0$ if and only if $\mathbb{P}_{\boldsymbol{\theta}_0} = \mathbb{P}_{\boldsymbol{\theta}}$.

If the distribution $\mathbb{P}$ does not belong to the parametric family $(\mathbb{P}_{\boldsymbol{\theta}}, \boldsymbol{\theta} \in \Theta)$, then the target of estimation can be naturally defined as the point of minimum of $-\mathbb{E} L(\boldsymbol{Y}, \boldsymbol{\theta})$. We will see that this point $\boldsymbol{\theta}_0$ indeed minimizes a special distance between the underlying measure $\mathbb{P}$ and the measures $\mathbb{P}_{\boldsymbol{\theta}}$ from the given parametric family.

The classical parametric statistical theory focuses mostly on asymptotic properties of the difference between $\widetilde{\boldsymbol{\theta}}$ and the true value $\boldsymbol{\theta}_0$ as the sample size $n$ tends to infinity. There is a vast literature on this issue. We only mention the book Ibragimov and Khas'minskij (1981), which provides a comprehensive study of asymptotic properties of maximum likelihood and Bayesian estimators. Typical results claim that the maximum likelihood and Bayes estimators are asymptotically optimal under certain regularity conditions. Large deviation results about minimum contrast estimators can be found in Jensen and Wood (1998) and Sieders and Dzhaparidze (1987), while subtle small sample size properties of these estimators are presented in Field (1982) and Field and Ronchetti (1990).

Another stream of the literature considers minimum contrast estimators in a general i.i.d. situation, when the parameter set $\Theta$ is a subset of some functional space. We mention the papers Van de Geer (1993), Birgé and Massart (1993), Birgé and Massart (1998), Birgé (2006) and references therein. The studies mostly focused on the concentration properties of the maximum $\max_{\boldsymbol{\theta}} L(\boldsymbol{Y}, \boldsymbol{\theta})$ rather on the properties of the estimator $\widetilde{\boldsymbol{\theta}}$ which is the point of maximum of $L(\boldsymbol{Y}, \boldsymbol{\theta})$. The established results are based on deep probabilistic facts from the empirical process theory; see e.g. van der Vaart and Wellner (1996). In this paper we also focus on the properties of the maximum of $L(\boldsymbol{Y}, \boldsymbol{\theta})$ over $\boldsymbol{\theta} \in \Theta$. However, we do not assume any particular structure of the contrast. Our basic result claims that if for every $\boldsymbol{\theta} \in \Theta$ the differences $L(\boldsymbol{Y}, \boldsymbol{\theta}) - L(\boldsymbol{Y}, \boldsymbol{\theta}_0)$ has exponential moments, then under rather general and mild conditions, the maximum $\max_{\boldsymbol{\theta}} \{L(\boldsymbol{Y}, \boldsymbol{\theta}) - L(\boldsymbol{Y}, \boldsymbol{\theta}_0)\}$ has similar exponential moments. In what follows, to keep notation shorter, we omit the argument $\boldsymbol{Y}$ in the contrast function $L(\boldsymbol{Y}, \boldsymbol{\theta})$ writing $L(\boldsymbol{\theta})$ instead of $L(\boldsymbol{Y}, \boldsymbol{\theta})$. However, one has to keep in mind that $L(\boldsymbol{\theta})$ is a random field that depends on the observed data $\boldsymbol{Y}$. We also denote

$$L(\boldsymbol{\theta}, \boldsymbol{\theta}_0) = L(\boldsymbol{\theta}) - L(\boldsymbol{\theta}_0).$$

To explain the main idea in this paper, introduce the function

$$\mathfrak{M}(\mu, \boldsymbol{\theta}, \boldsymbol{\theta}_0) \stackrel{\text{def}}{=} -\log \mathbb{E} \exp\{\mu L(\boldsymbol{\theta}, \boldsymbol{\theta}_0)\}.$$

Let $\mu^*$ be a maximizer of this function w.r.t. $\mu$, i.e.

$$\mu^*(\boldsymbol{\theta}) \stackrel{\text{def}}{=} \underset{\mu}{\operatorname{argmax}}\, \mathfrak{M}(\mu, \boldsymbol{\theta}, \boldsymbol{\theta}_0). \tag{1.3}$$

The *rate function* is defined via the Legendre transform of $L(\boldsymbol{\theta}, \boldsymbol{\theta}_0)$:

$$\mathfrak{M}^*(\boldsymbol{\theta}, \boldsymbol{\theta}_0) \stackrel{\text{def}}{=} \max_{\mu} \mathfrak{M}(\mu, \boldsymbol{\theta}, \boldsymbol{\theta}_0) = -\log \mathbb{E} \exp\{\mu^*(\boldsymbol{\theta}) L(\boldsymbol{\theta}, \boldsymbol{\theta}_0)\}. \tag{1.4}$$

Similar notions have already appeared in Chernoff (1952) and Bahadur (1960) for studying the models with i.i.d. observations.

Obviously $\mathfrak{M}^*(\boldsymbol{\theta}, \boldsymbol{\theta}_0) \geq \mathfrak{M}(0, \boldsymbol{\theta}, \boldsymbol{\theta}_0) = 0$. The following identity follows immediately from the above definition:

$$\mathbb{E} \exp\Big\{\mu^*(\boldsymbol{\theta}) L(\boldsymbol{\theta}, \boldsymbol{\theta}_0) + \mathfrak{M}^*(\boldsymbol{\theta}, \boldsymbol{\theta}_0)\Big\} = 1, \quad \boldsymbol{\theta} \in \Theta.$$

We aim to extend this pointwise identity to the supremum over $\boldsymbol{\theta} \in \Theta$, which particularly enables us to replace $\boldsymbol{\theta}$ with the estimator $\widetilde{\boldsymbol{\theta}}$. Unfortunately, in some situations, $\mathbb{E} \exp \sup_{\boldsymbol{\theta}} \{\mu^*(\boldsymbol{\theta}) L(\boldsymbol{\theta}, \boldsymbol{\theta}_0) + \mathfrak{M}^*(\boldsymbol{\theta}, \boldsymbol{\theta}_0)\} = \infty$. We illustrate this fact by some examples for a simple Gaussian liner model.

### 1.1. Examples for a linear Gaussian model

To illustrate how the quantities $\mu^*(\boldsymbol{\theta})$ and $\mathfrak{M}^*(\boldsymbol{\theta}, \boldsymbol{\theta}_0)$ can be computed let us consider the simplest case where $L(\boldsymbol{\theta}, \boldsymbol{\theta}_0)$ is a Gaussian field.

**Example 1.1.** [Gaussian contrast] Let for each pair $\boldsymbol{\theta}, \boldsymbol{\theta}' \in \Theta$, the difference $L(\boldsymbol{\theta}, \boldsymbol{\theta}') = L(\boldsymbol{\theta}) - L(\boldsymbol{\theta}')$ is a Gaussian random variable. In this case we call $L(\boldsymbol{\theta})$ a Gaussian contrast. With $M(\boldsymbol{\theta}, \boldsymbol{\theta}') = -\mathbb{E}L(\boldsymbol{\theta}, \boldsymbol{\theta}')$, $D^2(\boldsymbol{\theta}, \boldsymbol{\theta}') = \operatorname{Var} L(\boldsymbol{\theta}, \boldsymbol{\theta}')$, the random variable $L(\boldsymbol{\theta}, \boldsymbol{\theta}')$ is normal $\mathcal{N}(-M(\boldsymbol{\theta}, \boldsymbol{\theta}'), D^2(\boldsymbol{\theta}, \boldsymbol{\theta}'))$. Moreover,

$$\mathfrak{M}(\mu, \boldsymbol{\theta}, \boldsymbol{\theta}_0) = -\log \mathbb{E} \exp\{\mu L(\boldsymbol{\theta}, \boldsymbol{\theta}_0)\} = \mu M(\boldsymbol{\theta}, \boldsymbol{\theta}_0) - \mu^2 D^2(\boldsymbol{\theta}, \boldsymbol{\theta}_0)/2$$

and the values $\mu^*(\boldsymbol{\theta})$, $\mathfrak{M}^*(\boldsymbol{\theta}, \boldsymbol{\theta}_0)$ defined in (1.3)–(1.4) can be easily computed:

$$\mu^*(\boldsymbol{\theta}) = \operatorname*{argmax}_{\mu \geq 0}\{\mu M(\boldsymbol{\theta}, \boldsymbol{\theta}_0) - \mu^2 D^2(\boldsymbol{\theta}, \boldsymbol{\theta}_0)/2\} = \frac{M(\boldsymbol{\theta}, \boldsymbol{\theta}_0)}{D^2(\boldsymbol{\theta}, \boldsymbol{\theta}_0)},$$

$$\mathfrak{M}^*(\boldsymbol{\theta}, \boldsymbol{\theta}_0) = \sup_{\mu \geq 0} \mathfrak{M}(\mu, \boldsymbol{\theta}, \boldsymbol{\theta}_0) = \frac{M^2(\boldsymbol{\theta}, \boldsymbol{\theta}_0)}{2D^2(\boldsymbol{\theta}, \boldsymbol{\theta}_0)}.$$

The formula can be further simplified if $L(\boldsymbol{\theta})$ is a Gaussian log-likelihood.

**Example 1.2.** [Gaussian model] Let

$$L(\boldsymbol{\theta}, \boldsymbol{\theta}_0) = \log \frac{d\mathbb{P}_{\boldsymbol{\theta}}}{d\mathbb{P}_{\boldsymbol{\theta}_0}}(\boldsymbol{Y})$$

be a Gaussian random variable for any $\boldsymbol{\theta} \in \Theta$, and in addition $\mathbb{P} = \mathbb{P}_{\boldsymbol{\theta}_0}$ for some $\boldsymbol{\theta}_0 \in \Theta$. As in previous example, let $M(\boldsymbol{\theta}, \boldsymbol{\theta}_0)$ and $D(\boldsymbol{\theta}, \boldsymbol{\theta}_0)$ denote mean and variance of $L(\boldsymbol{\theta}, \boldsymbol{\theta}_0)$. The likelihood property implies $\mathbb{E}_{\boldsymbol{\theta}_0} \exp\{L(\boldsymbol{\theta}, \boldsymbol{\theta}_0)\} = 1$ yielding $M(\boldsymbol{\theta}, \boldsymbol{\theta}_0) = D^2(\boldsymbol{\theta}, \boldsymbol{\theta}_0)/2$ and hence, $\mu^*(\boldsymbol{\theta}) \equiv 1/2$ and $\mathfrak{M}^*(\boldsymbol{\theta}, \boldsymbol{\theta}_0) = M(\boldsymbol{\theta}, \boldsymbol{\theta}_0)/4$.

Finally we consider a classical linear Gaussian regression.

**Example 1.3.** [Linear Gaussian model] Consider the linear model $\boldsymbol{Y} = \boldsymbol{X}\boldsymbol{\theta}_0 + \sigma\boldsymbol{\varepsilon}$, where $\boldsymbol{Y} \in \mathbb{R}^n$, $\boldsymbol{\theta} \in \mathbb{R}^p$, $\boldsymbol{X}$ is a known $n \times p$ matrix, and $\boldsymbol{\varepsilon}$ is a white Gaussian noise in $\mathbb{R}^n$, i.e. $\varepsilon_i$ are i.i.d. standard normal. Then

$$L(\boldsymbol{\theta}) = -\|\boldsymbol{Y} - \boldsymbol{X}\boldsymbol{\theta}\|_n^2/(2\sigma^2),$$

where $\|\cdot\|_n$ denotes the standard Euclidian norm in $\mathbb{R}^n$. Obviously

$$M(\boldsymbol{\theta}, \boldsymbol{\theta}_0) = \|\boldsymbol{X}(\boldsymbol{\theta} - \boldsymbol{\theta}_0)\|_n^2/(2\sigma^2), \qquad D(\boldsymbol{\theta}, \boldsymbol{\theta}_0) = \|\boldsymbol{X}(\boldsymbol{\theta} - \boldsymbol{\theta}_0)\|_n^2/\sigma^2,$$

and thus (see Example 1.2)

$$\mathfrak{M}^*(\boldsymbol{\theta}, \boldsymbol{\theta}_0) = \|\boldsymbol{X}(\boldsymbol{\theta} - \boldsymbol{\theta}_0)\|_n^2/(8\sigma^2).$$

The log-likelihood ratio can be written as

$$L(\boldsymbol{\theta}, \boldsymbol{\theta}_0) = \langle \boldsymbol{X}(\boldsymbol{\theta} - \boldsymbol{\theta}_0), \boldsymbol{\varepsilon} \rangle_n/\sigma - \|\boldsymbol{X}(\boldsymbol{\theta} - \boldsymbol{\theta}_0)\|_n^2/(2\sigma^2).$$

Let $k$ denote the rank of the matrix $\boldsymbol{X}^\top \boldsymbol{X}$. Obviously $k \leq p$ and the vectors $\boldsymbol{X}(\boldsymbol{\theta} - \boldsymbol{\theta}_0)$ span a linear subspace $\mathcal{X}$ in $\mathbb{R}^n$ of dimension $k$. Denote by $\Pi$ the projector in $\mathbb{R}^n$ on $\mathcal{X}$. Then

$$\sup_{\boldsymbol{\theta} \in \mathbb{R}^p} \{\mu^*(\boldsymbol{\theta})L(\boldsymbol{\theta}, \boldsymbol{\theta}_0) + \mathfrak{M}^*(\boldsymbol{\theta}, \boldsymbol{\theta}_0)\}$$
$$= \sup_{\boldsymbol{\theta} \in \mathbb{R}^p} \left\{ \frac{\langle \boldsymbol{X}(\boldsymbol{\theta} - \boldsymbol{\theta}_0), \boldsymbol{\varepsilon} \rangle_n}{2\sigma} - \frac{\|\boldsymbol{X}(\boldsymbol{\theta} - \boldsymbol{\theta}_0)\|_n^2}{8\sigma^2} \right\}$$
$$= \sup_{\boldsymbol{u} \in \mathbb{R}^n} \left\{ \frac{\langle \Pi \boldsymbol{u}, \boldsymbol{\varepsilon} \rangle_n}{2\sigma} - \frac{\|\Pi \boldsymbol{u}\|_n^2}{8\sigma^2} \right\}$$
$$= \sup_{\boldsymbol{u} \in \mathbb{R}^n} \left\{ \frac{\langle \Pi \boldsymbol{u}, \Pi \boldsymbol{\varepsilon} \rangle_n}{2\sigma} - \frac{\|\Pi \boldsymbol{u}\|_n^2}{8\sigma^2} \right\} = \|\Pi \boldsymbol{\varepsilon}\|_n^2 / 2,$$

where the maximum is attained at any $\boldsymbol{u} \in \mathbb{R}^n$ such that $\Pi \boldsymbol{u} = 2\sigma \Pi \boldsymbol{\varepsilon}$. It is well known that $\|\Pi \boldsymbol{\varepsilon}\|_n^2$ follows $\chi^2$-distribution with $k$ degree of freedom and

$$\mathbb{E}_{\boldsymbol{\theta}_0} \exp \sup_{\boldsymbol{\theta}} \{\mu^*(\boldsymbol{\theta})L(\boldsymbol{\theta}, \boldsymbol{\theta}_0) + \mathfrak{M}^*(\boldsymbol{\theta}, \boldsymbol{\theta}_0)\} = \mathbb{E} \exp\{\|\Pi \boldsymbol{\varepsilon}\|_n^2/2\} = \infty.$$

However, for any positive $s < 1$, it holds by the same argument that

$$\sup_{\boldsymbol{\theta}} \{\mu^*(\boldsymbol{\theta})L(\boldsymbol{\theta}, \boldsymbol{\theta}_0) + s\mathfrak{M}^*(\boldsymbol{\theta}, \boldsymbol{\theta}_0)\}$$
$$= \sup_{\boldsymbol{u} \in \mathbb{R}^n} \left\{ \langle \Pi \boldsymbol{u}, \boldsymbol{\varepsilon} \rangle_n / (2\sigma) - (2-s)\|\Pi \boldsymbol{u}\|_k^2 / (8\sigma^2) \right\} = \|\Pi \boldsymbol{\varepsilon}\|_n^2 / (4 - 2s),$$

and thus

$$\mathbb{E}_{\boldsymbol{\theta}_0} \exp \sup_{\boldsymbol{\theta}} \{\mu^*(\boldsymbol{\theta})L(\boldsymbol{\theta}, \boldsymbol{\theta}_0) + s\mathfrak{M}^*(\boldsymbol{\theta}, \boldsymbol{\theta}_0)\} = \mathbb{E} \exp\left\{ \frac{\|\Pi \boldsymbol{\varepsilon}\|_n^2}{4 - 2s} \right\} = \left( \frac{2-s}{1-s} \right)^{k/2}.$$

An important feature of this inequality is that it only involves the effective dimension $k$ of the parameter space and does not depend on the design $\boldsymbol{X}$, noise level $\sigma^2$, sample size $n$, etc. Later we show that such a behaviour of the log-likelihood is not restricted to Gaussian linear models and it can be proved for a quite general statistical set-up.

### 1.2. Main result

The examples from Section 1.1 suggest to consider in the general situation the maximum of the random field $\mu^*(\boldsymbol{\theta})L(\boldsymbol{\theta}, \boldsymbol{\theta}_0) + s\mathfrak{M}^*(\boldsymbol{\theta}, \boldsymbol{\theta}_0)$ for $s < 1$. The main result of the paper shows that under some technical conditions this maximum is indeed stochastically bounded in a rather strong sense. Namely, for some $\rho \in (0, 1)$

$$\mathbb{E} \sup_{\boldsymbol{\theta} \in \Theta} \exp\left\{ \rho\left[\mu^*(\boldsymbol{\theta})L(\boldsymbol{\theta}, \boldsymbol{\theta}_0) + s\mathfrak{M}^*(\boldsymbol{\theta}, \boldsymbol{\theta}_0)\right] \right\} \leq C(\rho, s), \tag{1.5}$$

where $C(\rho, s)$ is a constant that can be easily controlled in typical examples. This result particularly yields that $\mu^*(\widetilde{\boldsymbol{\theta}})L(\widetilde{\boldsymbol{\theta}}, \boldsymbol{\theta}_0)$ and $\mathfrak{M}^*(\widetilde{\boldsymbol{\theta}}, \boldsymbol{\theta}_0)$ have bounded exponential moments. Another corollary of this fact is that $\widetilde{\boldsymbol{\theta}}$ concentrates on the sets $\mathcal{A}(\mathfrak{z}, \boldsymbol{\theta}_0) = \{\boldsymbol{\theta} : \mathfrak{M}^*(\boldsymbol{\theta}, \boldsymbol{\theta}_0) \leq \mathfrak{z}\}$ for sufficiently large $\mathfrak{z}$ in the sense that the probability $I\!\!P(\widetilde{\boldsymbol{\theta}} \notin \mathcal{A}(\mathfrak{z}, \boldsymbol{\theta}_0))$ is exponentially small in $\mathfrak{z}$. Usually every such concentration set is a root-n vicinity of the point $\boldsymbol{\theta}_0$. See Section 2.3 for precise formulations. Ibragimov and Khas'minskij (1981) stated a version of (1.5) for the i.i.d. case and used it to prove consistency of $\widetilde{\boldsymbol{\theta}}$.

We briefly comment on some useful features of the basic inequality (1.5). First of all this bound is non-asymptotic and may be used even if the sample size is small or moderate. It is also applicable in the situation when the parametric modeling assumption is misspecified. Our results may be used in such cases as well with the "true" parameter $\boldsymbol{\theta}_0$ defined as the maximum point of the contrast expected value: $\boldsymbol{\theta}_0 = \operatorname{argmax}_{\boldsymbol{\theta}} I\!\!E L(\boldsymbol{\theta})$.

Another interesting question is about the accuracy of estimation when the parameter set $\Theta$ is not compact. The typical results in the classical parametric theory has been established for compact parametric sets since this assumption simplifies considerably the conditions and the technical tools. There exist very few results for the case of non-compact sets. See Ibragimov and Khas'minskij (1981) for an example. Our conditions are quite mild and particularly, the parameter set can be non-compact and unbounded. Moreover, we present some examples in Section 4 illustrating that the quality of the minimum contrast estimation can heavily depend on topological properties of $\Theta$ and on the behavior of the rate function $\mathfrak{M}^*(\boldsymbol{\theta}, \boldsymbol{\theta}_0)$ for large $\boldsymbol{\theta}$. The corresponding accuracy of estimation can be different from the classical root-$n$ behavior.

The paper is organized as follows. The main result is presented in Section 2. Section 2.3 presents some useful corollaries of (1.5) describing concentration properties of $\widetilde{\boldsymbol{\theta}}$, some risk bounds, confidence sets for the target parameter $\boldsymbol{\theta}_0$ based on the $L(\widetilde{\boldsymbol{\theta}}, \boldsymbol{\theta})$. Section 2.4 specifies the approach to the important case of a smooth contrast. In this situation the main conditions ensuring (1.5) are substantially simplified. Section 3 illustrates how our approach applies to the classical i.i.d. case while Section 4 presents some applications of the general exponential bound to three particular problems: estimation of the median, of the scale parameter of an exponential model and of the change point location. Although these examples have already been studied, the proposed approach reveals some new features of the classical least squares and least absolute deviation estimators in the cases when the parametric assumption is misspecified or the parameter set is not compact. In the case of median estimation the result applies even if the observations do not have the first moment. The last example in this section considers the prominent change point problem. We particularly show that in the case when the size of the jump is completely unknown, the accuracy of estimation of its location differs from the well known parametric rate $1/n$ and it depends on the distance of the change point to the edge of the observation interval and involves an extra iterated-log factor.

## 2. Risk bound for the minimum contrast

This section presents a general exponential bound on the minimum contrast value in a rather general set-up. Let $-L(\boldsymbol{\theta}),\ \boldsymbol{\theta}\in\Theta$, be a random *contrast function* of a finite dimensional parameter $\boldsymbol{\theta}\in\Theta\subset I\!\!R^p$ given on some probability space $(\Omega,\mathcal{F},I\!\!P)$. We also assume that $L(\boldsymbol{\theta})$ is separable random field and $I\!\!E L(\boldsymbol{\theta})$ exists for all $\boldsymbol{\theta}\in\Theta$. The *minimum contrast estimator* is defined as a minimizer of $-L(\boldsymbol{\theta})$ and the target of estimation is the value $\boldsymbol{\theta}_0$ which minimizes the expectation $-I\!\!E L(\boldsymbol{\theta})$. It is clear that for any $\boldsymbol{\theta}^\circ\in\Theta$

$$\widetilde{\boldsymbol{\theta}}=\operatorname*{argmax}_{\boldsymbol{\theta}\in\Theta}L(\boldsymbol{\theta},\boldsymbol{\theta}^\circ)\quad\text{and}\quad\boldsymbol{\theta}_0=\operatorname*{argmax}_{\boldsymbol{\theta}\in\Theta}I\!\!E L(\boldsymbol{\theta},\boldsymbol{\theta}^\circ).$$

Our study focuses on the value of maximum in $\boldsymbol{\theta}$ of the random field $L(\boldsymbol{\theta},\boldsymbol{\theta}_0)$:

$$L(\widetilde{\boldsymbol{\theta}},\boldsymbol{\theta}_0)=\sup_{\boldsymbol{\theta}\in\Theta}L(\boldsymbol{\theta},\boldsymbol{\theta}_0)=\sup_{\boldsymbol{\theta}\in\Theta}\{L(\boldsymbol{\theta})-L(\boldsymbol{\theta}_0)\}.$$

By definition, $L(\widetilde{\boldsymbol{\theta}},\boldsymbol{\theta}_0)$ is a non-negative random variable.

### 2.1. Preliminaries. The case of a discrete parameter set

The main goal of this paper is to obtain exponential bounds for the supremum in $\boldsymbol{\theta}$ of the random field $L(\boldsymbol{\theta},\boldsymbol{\theta}_0)$, without specifying a particular structure of the model or contrast function $L(\boldsymbol{\theta})$. Instead we impose some conditions of finite exponential moments for the increments $L(\boldsymbol{\theta},\boldsymbol{\theta}')=L(\boldsymbol{\theta})-L(\boldsymbol{\theta}')$. With $\mathfrak{M}(\mu,\boldsymbol{\theta},\boldsymbol{\theta}_0)=-\log I\!\!E\exp\{\mu L(\boldsymbol{\theta},\boldsymbol{\theta}_0)\}$, the *global exponential moment* condition reads as follows:

**(EG)** *For any $\boldsymbol{\theta}\in\Theta$ the set $\Upsilon(\boldsymbol{\theta})=\{\mu\in(0,\infty):\ \mathfrak{M}(\mu,\boldsymbol{\theta},\boldsymbol{\theta}_0)<\infty\}$ is non-empty.*

Note that $\Upsilon(\boldsymbol{\theta})$ is an interval because $\mathfrak{M}(\mu,\boldsymbol{\theta},\boldsymbol{\theta}_0)<\infty$ implies $\mathfrak{M}(\mu',\boldsymbol{\theta},\boldsymbol{\theta}_0)<\infty$ for all $\mu'<\mu$. Moreover, in the basic example of the log-likelihood contrast, it holds $\mathfrak{M}(1,\boldsymbol{\theta},\boldsymbol{\theta}_0)=-\log I\!\!E_{\boldsymbol{\theta}_0}(dI\!\!P_{\boldsymbol{\theta}}/dI\!\!P_{\boldsymbol{\theta}_0})\leq 0$ for all $\boldsymbol{\theta}$ and the condition $(EG)$ is fulfilled automatically with $(0,1]\subset\Upsilon(\boldsymbol{\theta},\boldsymbol{\theta}_0)$.

Under the condition $(EG)$ the functions $\mu^*(\boldsymbol{\theta})$ and $\mathfrak{M}^*(\boldsymbol{\theta},\boldsymbol{\theta}_0)$ from (1.3)–(1.4) are non-trivial and correctly defined. Usually these functions can be easily evaluated in a small neighborhood of the target parameter $\boldsymbol{\theta}_0$. However, it might be difficult to compute them for all $\boldsymbol{\theta}\in\Theta$. Therefore, in the sequel we proceed with another function $\mu(\boldsymbol{\theta})$, which can be viewed as a rough approximation of $\mu^*(\boldsymbol{\theta})$. Section 4 provides some examples. So, let $\mu(\boldsymbol{\theta})$ be a given function taking values in $\Upsilon(\boldsymbol{\theta},\boldsymbol{\theta}_0)$. Define

$$\mathfrak{M}(\boldsymbol{\theta},\boldsymbol{\theta}_0)\stackrel{\text{def}}{=}\mathfrak{M}(\mu(\boldsymbol{\theta}),\boldsymbol{\theta},\boldsymbol{\theta}_0)=-\log I\!\!E\exp\{\mu(\boldsymbol{\theta})L(\boldsymbol{\theta},\boldsymbol{\theta}_0)\}.$$

The most important requirement on $\mu(\boldsymbol{\theta})$ is that $\mathfrak{M}(\boldsymbol{\theta},\boldsymbol{\theta}_0)$ is positive and increases as $\boldsymbol{\theta}$ moves away from $\boldsymbol{\theta}_0$. By definition, for any $\boldsymbol{\theta}\in\Theta$,

$$I\!\!E\exp\{\mu(\boldsymbol{\theta})L(\boldsymbol{\theta},\boldsymbol{\theta}_0)+\mathfrak{M}(\boldsymbol{\theta},\boldsymbol{\theta}_0)\}=1. \tag{2.1}$$

This means that the random function $\mu(\boldsymbol{\theta})L(\boldsymbol{\theta},\boldsymbol{\theta}_0) + \mathfrak{M}(\boldsymbol{\theta},\boldsymbol{\theta}_0)$ has bounded exponential moments for every $\boldsymbol{\theta}$. We aim to derive a similar fact for the supremum of this function in $\boldsymbol{\theta} \in \Theta$. More precisely, we are interested in bounding the following value:

$$\mathfrak{Q}(\rho,s) \stackrel{\text{def}}{=} I\!\!E \sup_{\boldsymbol{\theta} \in \Theta} \exp\Big\{\rho\big[\mu(\boldsymbol{\theta})L(\boldsymbol{\theta},\boldsymbol{\theta}_0) + s\mathfrak{M}(\boldsymbol{\theta},\boldsymbol{\theta}_0)\big]\Big\}, \qquad (2.2)$$

where $\rho, s \in [0,1]$.

We begin with a rough upper bound for a special case of a discrete parameter set.

**Theorem 2.1.** *Assume* (EG) *and let* $\Theta$ *be a discrete set. Then for any* $s < 1$

$$\begin{aligned}\mathfrak{Q}(1,s) &= I\!\!E \sup_{\boldsymbol{\theta} \in \Theta} \exp\big\{\mu(\boldsymbol{\theta})L(\boldsymbol{\theta},\boldsymbol{\theta}_0) + s\mathfrak{M}(\boldsymbol{\theta},\boldsymbol{\theta}_0)\big\} \\ &\leq \sum_{\boldsymbol{\theta} \in \Theta} \exp\big\{-(1-s)\mathfrak{M}(\boldsymbol{\theta},\boldsymbol{\theta}_0)\big\}. \end{aligned} \qquad (2.3)$$

*Proof.* Since $I\!\!E \exp\big\{\mu(\boldsymbol{\theta})L(\boldsymbol{\theta},\boldsymbol{\theta}_0) + s\mathfrak{M}(\boldsymbol{\theta},\boldsymbol{\theta}_0)\big\} = \exp\big\{-(1-s)\mathfrak{M}(\boldsymbol{\theta},\boldsymbol{\theta}_0)\big\}$, we obviously have

$$\mathfrak{Q}(1,s) \leq \sum_{\boldsymbol{\theta} \in \Theta} I\!\!E \exp\big\{\mu(\boldsymbol{\theta})L(\boldsymbol{\theta},\boldsymbol{\theta}_0) + s\mathfrak{M}(\boldsymbol{\theta},\boldsymbol{\theta}_0)\big\} = \sum_{\boldsymbol{\theta} \in \Theta} \exp\big\{-(1-s)\mathfrak{M}(\boldsymbol{\theta},\boldsymbol{\theta}_0)\big\}.$$

$\square$

Usually, the function $\mathfrak{M}(\boldsymbol{\theta},\boldsymbol{\theta}_0)$ rapidly grows as $\boldsymbol{\theta}$ moves away from $\boldsymbol{\theta}_0$. This property is often sufficient to bound the sum in the right hand-side of (2.3) by a fixed constant.

Although Theorem 2.1 is a rather simple corollary of (2.1), the bound (2.3) yields a number of useful statistical corollaries. Some of them are presented in Section 2.3. However, even in discrete case, this bound may be too rough (see the example in Section 4.3). It is also clear that (2.3) is useless in the continuous case. The next section demonstrates how the bound (2.3) can be extended to the case of an arbitrary parameter set.

### 2.2. The general exponential bound

Here we aim to extend the exponential bound (2.3) from the discrete case to the case of an arbitrary finite dimensional parameter set. We apply the standard approach which evaluates the supremum over the whole parameter set $\Theta$ via a weighted sum of local maxima.

Define for any $\boldsymbol{\theta}, \boldsymbol{\theta}' \in \Theta$

$$\zeta(\boldsymbol{\theta}) \stackrel{\text{def}}{=} \mu(\boldsymbol{\theta})\big\{L(\boldsymbol{\theta},\boldsymbol{\theta}_0) - I\!\!E L(\boldsymbol{\theta},\boldsymbol{\theta}_0)\big\}, \qquad \zeta(\boldsymbol{\theta},\boldsymbol{\theta}') \stackrel{\text{def}}{=} \zeta(\boldsymbol{\theta}) - \zeta(\boldsymbol{\theta}').$$

Note that the dependence of $\zeta(\boldsymbol{\theta},\boldsymbol{\theta}')$ on $\boldsymbol{\theta}_0$ disappears if $\mu(\boldsymbol{\theta}) = \mu(\boldsymbol{\theta}')$.

Usually the local properties of the centered contrast difference $\zeta(\boldsymbol{\theta}, \boldsymbol{\theta}')$ are controlled by the variance $D^2(\boldsymbol{\theta}, \boldsymbol{\theta}') = \operatorname{Var} \zeta(\boldsymbol{\theta}, \boldsymbol{\theta}')$, which defines a semi-metric on $\Theta$ see, e.g. van der Vaart and Wellner (1996). However, in some cases, it is more convenient to deal with a slightly different metric which we denote by $\mathfrak{S}(\boldsymbol{\theta}, \boldsymbol{\theta}')$. This metric usually bounds the standard deviation $D(\boldsymbol{\theta}, \boldsymbol{\theta}')$ from above. Sections 2.4 and 3 present some typical examples of constructing such a metric. Below in this section we assume that the metric $\mathfrak{S}(\cdot, \cdot)$ is given. Define for any point $\boldsymbol{\theta}^\circ \in \Theta$ and a radius $\epsilon > 0$ the ball

$$\mathcal{B}(\epsilon, \boldsymbol{\theta}^\circ) = \{\boldsymbol{\theta} : \mathfrak{S}(\boldsymbol{\theta}, \boldsymbol{\theta}^\circ) \leq \epsilon\}.$$

To control the local behavior of the process $L(\boldsymbol{\theta})$ within any such ball $\mathcal{B}(\epsilon, \boldsymbol{\theta}^\circ)$, we impose the following *local exponential* condition:

**(EL)** There exist $\epsilon > 0$ and $\overline{\lambda} > 0$ such that for any $\boldsymbol{\theta}^\circ \in \Theta$, $\nu_0 > 0$, and $\lambda \leq \overline{\lambda}$

$$\sup_{\boldsymbol{\theta}, \boldsymbol{\theta}' \in \mathcal{B}(\epsilon, \boldsymbol{\theta}^\circ)} \log I\!E \exp\{2\lambda \xi(\boldsymbol{\theta}, \boldsymbol{\theta}')\} \leq 2\nu_0^2 \lambda^2,$$

where

$$\xi(\boldsymbol{\theta}, \boldsymbol{\theta}') \stackrel{\text{def}}{=} \frac{\zeta(\boldsymbol{\theta}, \boldsymbol{\theta}')}{\mathfrak{S}(\boldsymbol{\theta}, \boldsymbol{\theta}')}.$$

In fact, this condition only requires that every random increment $\xi(\boldsymbol{\theta}, \boldsymbol{\theta}')$ has bounded exponential moment for some $\overline{\lambda} > 0$. Then Lemma 5.8 from the Appendix implies the prescribed quadratic behavior in $\lambda$ for $\lambda \leq \overline{\lambda}$.

For a fixed $\boldsymbol{\theta}^\circ \in \Theta$ and $\epsilon' \leq \epsilon$, by $\mathbb{N}(\epsilon', \epsilon, \boldsymbol{\theta}^\circ)$ we denote the local covering number defined as the minimal number of balls $\mathcal{B}(\epsilon', \cdot)$ required to cover the ball $\mathcal{B}(\epsilon, \boldsymbol{\theta}^\circ)$. With this covering number we associate the *local entropy*

$$\mathbb{Q}(\epsilon, \boldsymbol{\theta}^\circ) \stackrel{\text{def}}{=} \sum_{k=1}^{\infty} 2^{-k} \log \mathbb{N}(2^{-k}\epsilon, \epsilon, \boldsymbol{\theta}^\circ).$$

We begin with a local result which bounds the maximum of the process $L(\boldsymbol{\theta})$ over a local ball $\mathcal{B}(\epsilon, \boldsymbol{\theta}^\circ)$.

**Theorem 2.2.** *Assume* (EG) *and* (EL) *with some* $\epsilon > 0$, $\nu_0 \geq 0$, *and* $\overline{\lambda} > 0$. *Let also* $\rho < 1$ *be such that* $\rho\epsilon/(1-\rho) \leq \overline{\lambda}$. *Then for any* $\boldsymbol{\theta}^\circ \in \Theta$

$$\log I\!E \sup_{\boldsymbol{\theta} \in \mathcal{B}(\epsilon, \boldsymbol{\theta}^\circ)} \exp\left\{\rho\bigl[\mu(\boldsymbol{\theta})L(\boldsymbol{\theta}, \boldsymbol{\theta}_0) + \mathfrak{M}(\boldsymbol{\theta}, \boldsymbol{\theta}_0)\bigr]\right\} \leq \frac{2\nu_0^2 \epsilon^2 \rho^2}{1-\rho} + (1-\rho)\mathbb{Q}(\epsilon, \boldsymbol{\theta}^\circ).$$

The next theorem is the global bound which generalizes the upper bound from Theorem 2.1.

**Theorem 2.3.** *Assume* (EG) *and* (EL) *for some* $\overline{\lambda}, \nu_0, \epsilon$, *and let* $\pi(\cdot)$ *be a $\sigma$-finite measure on* $\Theta$ *such that*

$$\sup_{\boldsymbol{\theta} \in \mathcal{B}(\epsilon, \boldsymbol{\theta}^\circ)} \frac{\pi\bigl(\mathcal{B}(\epsilon, \boldsymbol{\theta})\bigr)}{\pi\bigl(\mathcal{B}(\epsilon, \boldsymbol{\theta}^\circ)\bigr)} \leq \nu_1 \qquad (2.4)$$

for some $\nu_1 \in [1, \infty)$. Let for some $\rho, s < 1$, it holds $\rho\epsilon/(1-\rho) \leq \overline{\lambda}$ and the function $\mathfrak{M}_\epsilon(\boldsymbol{\theta}^\circ, \boldsymbol{\theta}_0) = \inf_{\boldsymbol{\theta} \in \mathcal{B}(\epsilon, \boldsymbol{\theta}^\circ)} \mathfrak{M}(\boldsymbol{\theta}, \boldsymbol{\theta}_0)$ fulfill

$$\mathfrak{H}_\epsilon(\rho, s) \stackrel{\text{def}}{=} \log \left( \int_\Theta \frac{1}{\pi(\mathcal{B}(\epsilon, \boldsymbol{\theta}))} \exp\{-\rho(1-s)\mathfrak{M}_\epsilon(\boldsymbol{\theta}, \boldsymbol{\theta}_0)\} \pi(d\boldsymbol{\theta}) \right) < \infty. \quad (2.5)$$

Let finally $\mathbb{Q}(\epsilon, \boldsymbol{\theta}^\circ) \leq \overline{\mathbb{Q}}(\epsilon)$ for all $\boldsymbol{\theta} \in \Theta$. Then the value $\mathfrak{Q}(\rho, s)$ from (2.2) satisfies

$$\mathfrak{Q}(\rho, s) \leq \frac{2\nu_0^2 \epsilon^2 \rho^2}{1 - \rho} + (1-\rho)\overline{\mathbb{Q}}(\epsilon) + \log(\nu_1) + \mathfrak{H}_\epsilon(\rho, s). \quad (2.6)$$

As in Theorem 2.1, proper growth conditions on the function $\mathfrak{M}(\boldsymbol{\theta}, \boldsymbol{\theta}_0)$ ensure that the integral $\mathfrak{H}_\epsilon(\rho, s)$ in (2.6) is bounded by a fixed constant.

### 2.3. Some corollaries

This section demonstrates how Theorems 2.1–2.3 can be used in the statistical analysis of the minimum contrast estimator $\widetilde{\boldsymbol{\theta}} = \mathrm{argmax}_{\boldsymbol{\theta} \in \Theta} L(\boldsymbol{\theta})$. We show that probabilistic properties of this estimator may be easily derived from the following inequality: for prescribed $\rho, s < 1$,

$$I\!\!E \exp\left\{ \rho\big[\mu(\widetilde{\boldsymbol{\theta}}) L(\widetilde{\boldsymbol{\theta}}, \boldsymbol{\theta}_0) + s\mathfrak{M}(\widetilde{\boldsymbol{\theta}}, \boldsymbol{\theta}_0)\big] \right\} \leq \mathfrak{Q}(\rho, s), \quad (2.7)$$

which obviously follows from Theorem 2.3 and the definition (2.2) of $\mathfrak{Q}(\rho, s)$.

#### 2.3.1. A risk bound for the "natural" loss

A first corollary of Theorem 2.1 presents exponential bounds separately for the minimum contrast $L(\widetilde{\boldsymbol{\theta}}, \boldsymbol{\theta}_0)$ and for the "natural" loss $\mathfrak{M}(\widetilde{\boldsymbol{\theta}}, \boldsymbol{\theta}_0)$.

**Corollary 2.4.** *For any $\rho, s < 1$*

$$I\!\!E \exp\left\{ \rho\mu(\widetilde{\boldsymbol{\theta}}) L(\widetilde{\boldsymbol{\theta}}, \boldsymbol{\theta}_0) \right\} \leq \mathfrak{Q}(\rho, 0), \quad (2.8)$$

$$I\!\!E \exp\left\{ \rho s\, \mathfrak{M}(\widetilde{\boldsymbol{\theta}}, \boldsymbol{\theta}_0) \right\} \leq \mathfrak{Q}(\rho, s). \quad (2.9)$$

Substituting $s = 0$ in (2.7) yields the first bound. To prove the second one, notice that $L(\widetilde{\boldsymbol{\theta}}, \boldsymbol{\theta}_0) \geq 0$. Therefore the elementary inequality $\mathbf{1}\{x \geq 0\} \leq \exp(\mu x)$ for any $\mu > 0$ yields (see also (2.7))

$$\begin{aligned}
I\!\!E \exp\{\rho s\, \mathfrak{M}(\widetilde{\boldsymbol{\theta}}, \boldsymbol{\theta}_0)\} &= I\!\!E \exp\{\rho s\, \mathfrak{M}(\widetilde{\boldsymbol{\theta}}, \boldsymbol{\theta}_0)\} \mathbf{1}\{L(\widetilde{\boldsymbol{\theta}}, \boldsymbol{\theta}_0) \geq 0\} \\
&\leq I\!\!E \exp\{\rho s\, \mathfrak{M}(\widetilde{\boldsymbol{\theta}}, \boldsymbol{\theta}_0) + \rho\mu(\widetilde{\boldsymbol{\theta}}) L(\widetilde{\boldsymbol{\theta}}, \boldsymbol{\theta}_0)\} \leq \mathfrak{Q}(\rho, s).
\end{aligned}$$

Notice the exponential bound (2.9) implies a similar risk bound for a polynomial loss $\left|\mathfrak{M}(\widetilde{\boldsymbol{\theta}}, \boldsymbol{\theta}_0)\right|^r$; see Lemma 5.7 for a precise result.

### 2.3.2. Concentration properties of the estimator $\widetilde{\boldsymbol{\theta}}$

The assertion (2.7) can be used for establishing the concentration property of the estimator $\widetilde{\boldsymbol{\theta}}$. Consider the sets

$$\mathcal{A}(r, \boldsymbol{\theta}_0) \stackrel{\text{def}}{=} \{\boldsymbol{\theta} : \mathfrak{M}(\boldsymbol{\theta}, \boldsymbol{\theta}_0) \leq r\}$$

for some $r > 0$. The next result shows that the estimator $\widetilde{\boldsymbol{\theta}}$ leaves the set $\mathcal{A}(r, \boldsymbol{\theta}_0)$ with the exponentially small probability of order $\exp(-\rho s r)$.

**Corollary 2.5.** *For any $\rho, s < 1$, it holds*

$$I\!P\big(\widetilde{\boldsymbol{\theta}} \notin \mathcal{A}(r, \boldsymbol{\theta}_0)\big) \leq \mathfrak{Q}(\rho, s) \exp(-\rho s r).$$

*Proof.* The inequalities $L(\widetilde{\boldsymbol{\theta}}, \boldsymbol{\theta}_0) \geq 0$ and $\mathfrak{M}(\widetilde{\boldsymbol{\theta}}, \boldsymbol{\theta}_0) > r$ for $\widetilde{\boldsymbol{\theta}} \notin \mathcal{A}(r, \boldsymbol{\theta}_0)$ imply

$$I\!E e^{\rho s r} \mathbf{1}\big(\widetilde{\boldsymbol{\theta}} \notin \mathcal{A}(r, \boldsymbol{\theta}_0)\big) \leq I\!E \exp\Big\{\rho\big[\mu(\widetilde{\boldsymbol{\theta}}) L(\widetilde{\boldsymbol{\theta}}, \boldsymbol{\theta}_0) + s \mathfrak{M}(\widetilde{\boldsymbol{\theta}}, \boldsymbol{\theta}_0)\big]\Big\} \leq \mathfrak{Q}(\rho, s)$$

and the assertion follows. □

In typical situations, $\mathfrak{M}(\boldsymbol{\theta}, \boldsymbol{\theta}_0)$ is proportional to the sample size $n$ and each set $\mathcal{A}(r, \boldsymbol{\theta}_0)$ corresponds to a root-$n$ neighborhood of the point $\boldsymbol{\theta}_0$. See the Section 3 for applications related to the i.i.d. case.

### 2.3.3. Confidence sets based on $L(\widetilde{\boldsymbol{\theta}}, \boldsymbol{\theta})$

Next we discuss how the exponential bound (2.7) can be used for constructing the confidence sets for the target $\boldsymbol{\theta}_0$ based on the optimized contrast $L(\widetilde{\boldsymbol{\theta}}, \boldsymbol{\theta})$. The inequality (2.8) claims that $L(\widetilde{\boldsymbol{\theta}}, \boldsymbol{\theta}_0)$ is stochastically bounded. This justifies the following construction of confidence sets:

$$\mathcal{E}(\mathfrak{z}) = \big\{\boldsymbol{\theta} \in \Theta : L(\widetilde{\boldsymbol{\theta}}, \boldsymbol{\theta}) \leq \mathfrak{z}\big\}.$$

To evaluate the covering probability, consider first the case when $\mu(\boldsymbol{\theta}) \geq \mu_* > 0$ uniformly in $\boldsymbol{\theta} \in \Theta$. The next result claims that $\mathcal{E}(\mathfrak{z})$ does not cover the true value $\boldsymbol{\theta}_0$ with a probability which decreases exponentially with $\mathfrak{z}$.

**Corollary 2.6.** *Assume that $\mu(\boldsymbol{\theta}) \geq \mu_* > 0$. Then for any $\mathfrak{z} > 0$ and any $\rho < 1$*

$$I\!P\big(\boldsymbol{\theta}_0 \notin \mathcal{E}(\mathfrak{z})\big) \leq \mathfrak{Q}(\rho, 0) \exp\big\{-\rho \mu_* \mathfrak{z}\big\}.$$

*Proof.* The bound (2.8) implies

$$\begin{aligned}
I\!P\big(\boldsymbol{\theta}_0 \notin \mathcal{E}(\mathfrak{z})\big) &= I\!P\big(L(\widetilde{\boldsymbol{\theta}}, \boldsymbol{\theta}_0) > \mathfrak{z}\big) \\
&\leq I\!E \exp\big\{-\rho \mu(\widetilde{\boldsymbol{\theta}}) \mathfrak{z}\big\} \exp\big\{\rho \mu(\widetilde{\boldsymbol{\theta}}) L(\widetilde{\boldsymbol{\theta}}, \boldsymbol{\theta}_0)\big\} \\
&\leq \exp\big\{-\rho \mu_* \mathfrak{z}\big\} I\!E \exp\big\{\rho \mu(\widetilde{\boldsymbol{\theta}}) L(\widetilde{\boldsymbol{\theta}}, \boldsymbol{\theta}_0)\big\} \\
&\leq \mathfrak{Q}(\rho, 0) \exp\big\{-\rho \mu_* \mathfrak{z}\big\}
\end{aligned}$$

as required. □

In the case when the function $\mu(\boldsymbol{\theta})$ cannot be uniformly bounded from below by a positive constant, we assume that such a bound exists for every set $\mathcal{A}(r, \boldsymbol{\theta}_0)$. Denote

$$\mu_*(r) \stackrel{\text{def}}{=} \inf_{\boldsymbol{\theta} \in \mathcal{A}(r, \boldsymbol{\theta}_0)} \mu(\boldsymbol{\theta}).$$

Then

$$I\!P\big(\boldsymbol{\theta}_0 \notin \mathcal{E}(\mathfrak{z})\big) \leq I\!P\big(\boldsymbol{\theta}_0 \notin \mathcal{E}(\mathfrak{z}), \widetilde{\boldsymbol{\theta}} \in \mathcal{A}(r, \boldsymbol{\theta}_0)\big) + I\!P\big(\widetilde{\boldsymbol{\theta}} \notin \mathcal{A}(r, \boldsymbol{\theta}_0)\big)$$

and combining Corollaries 2.5–2.6 yields

**Corollary 2.7.** *For any $\mathfrak{z} > 0$ and any $\rho, s < 1$ and any $r > 0$*

$$I\!P\big(\boldsymbol{\theta}_0 \notin \mathcal{E}(\mathfrak{z})\big) \leq \mathfrak{Q}(\rho, 0) \exp\{-\rho\mu_*(r)\mathfrak{z}\} + \mathfrak{Q}(\rho, s) \exp\{-\rho s r\}.$$

A reasonable choice of $r$ in this bound is given by the balance relation $\mu_*(r)\mathfrak{z} = sr$. With this choice the bound of Corollary 2.6 may by replaced by

$$I\!P\big(\boldsymbol{\theta}_0 \notin \mathcal{E}(\mathfrak{z})\big) \leq 2\mathfrak{Q}(\rho, s) \exp\{-\rho\mu_*(r)\mathfrak{z}\}.$$

## 2.4. Exponential bounds for smooth contrasts

This section deals with the case when the contrast $L(\boldsymbol{\theta})$ is a smooth function of $\boldsymbol{\theta}$. In this situation, the local condition $(EL)$ is easy to verify. Moreover, the local balls $\mathcal{B}(\epsilon, \boldsymbol{\theta})$ nearly coincide with usual Euclidean ellipsoids and the local entropy can be easily bounded by an absolute constant only depending on the dimensionality $p$ of the parameter space $\Theta$.

Suppose $\Theta$ is a convex set in $I\!R^p$ and the function $L(\boldsymbol{\theta})$ along with the scaling factor $\mu(\boldsymbol{\theta})$ are differentiable w.r.t. $\boldsymbol{\theta}$. Below, the symbol $\nabla$ stands for the gradient w.r.t. $\boldsymbol{\theta}$.

Define

$$V(\boldsymbol{\theta}) \stackrel{\text{def}}{=} I\!E \nabla \zeta(\boldsymbol{\theta}) \big[\nabla \zeta(\boldsymbol{\theta})\big]^\top,$$

$$H(\lambda, \gamma, \boldsymbol{\theta}) \stackrel{\text{def}}{=} \log I\!E \exp\left\{2\lambda \frac{\gamma^\top \nabla \zeta(\boldsymbol{\theta})}{\sqrt{\gamma^\top V(\boldsymbol{\theta})\gamma}}\right\}.$$

for every unit vector $\gamma \in I\!R^p$. To simplify the presentation, here and in what follows we assume that every matrix $V(\boldsymbol{\theta})$ is non-degenerated. It is easy to see that $H(0, \gamma, \boldsymbol{\theta}) = 0$, $\partial H(0, \gamma, \boldsymbol{\theta})/\partial \lambda = 0$, and

$$\left.\frac{\partial^2 H(\lambda, \gamma, \boldsymbol{\theta})}{\partial^2 \lambda}\right|_{\lambda=0} = \frac{4\gamma^\top I\!E \nabla \zeta(\boldsymbol{\theta})\big[\nabla \zeta(\boldsymbol{\theta})\big]^\top \gamma}{\gamma^\top V(\boldsymbol{\theta})\gamma} = 4.$$

Therefore for small $\lambda$ $H(\lambda, \gamma, \boldsymbol{\theta}) \approx 2\lambda^2$. Below we assume that this property is fulfilled uniformly in $\boldsymbol{\theta} \in \Theta$ and in $\gamma$ over the unit sphere $S^p$ in $I\!R^p$.

**(ED)** There exists $\overline{\lambda} > 0$ such that for some $\nu_0 \geq 1$ uniformly in $\boldsymbol{\theta} \in \Theta$

$$\sup_{|\lambda| \leq \overline{\lambda}} \sup_{\gamma \in S^p} \lambda^{-2} H(\lambda, \gamma, \boldsymbol{\theta}) \leq 2\nu_0^2. \tag{2.10}$$

Now we define the metric $\mathfrak{S}(\boldsymbol{\theta}, \boldsymbol{\theta}')$ by

$$\mathfrak{S}^2(\boldsymbol{\theta}, \boldsymbol{\theta}') \stackrel{\text{def}}{=} \sup_{t \in [0,1]} (\boldsymbol{\theta} - \boldsymbol{\theta}')^\top V\big[(1-t)\boldsymbol{\theta}' + t\boldsymbol{\theta}\big](\boldsymbol{\theta} - \boldsymbol{\theta}'). \tag{2.11}$$

Define also for every $\boldsymbol{\theta}^\circ \in \Theta$ and $\epsilon > 0$ the ellipsoid $\mathcal{B}'(\epsilon, \boldsymbol{\theta}^\circ)$ by

$$\mathcal{B}'(\epsilon, \boldsymbol{\theta}^\circ) = \Big\{\boldsymbol{\theta} : (\boldsymbol{\theta} - \boldsymbol{\theta}^\circ)^\top V(\boldsymbol{\theta}^\circ)(\boldsymbol{\theta} - \boldsymbol{\theta}^\circ) \leq \epsilon^2\Big\}.$$

Obviously $\mathcal{B}(\epsilon, \boldsymbol{\theta}^\circ) \subseteq \mathcal{B}'(\epsilon, \boldsymbol{\theta}^\circ)$.

In what follows, we assume that the radius $\epsilon$ can be chosen in such a way that the functions $V(\boldsymbol{\theta})$ and $\mathfrak{M}(\boldsymbol{\theta}, \boldsymbol{\theta}_0)$ have bounded fluctuations within the ball $\mathcal{B}'(\epsilon, \boldsymbol{\theta}^\circ)$ for every $\boldsymbol{\theta}^\circ \in \Theta$. More precisely, for a given function $f(\cdot)$ define its magnitude over $\mathcal{B}'(\epsilon, \boldsymbol{\theta}^\circ)$ by

$$\mathfrak{A}_\epsilon f(\boldsymbol{\theta}^\circ) \stackrel{\text{def}}{=} \sup_{\boldsymbol{\theta}, \boldsymbol{\theta}' \in \mathcal{B}'(\epsilon, \boldsymbol{\theta}^\circ)} \frac{f(\boldsymbol{\theta})}{f(\boldsymbol{\theta}')}.$$

Similarly, the magnitude of the matrix $V(\boldsymbol{\theta})$ over $\mathcal{B}'(\epsilon, \boldsymbol{\theta}^\circ)$ is computed as follows

$$\mathfrak{A}_\epsilon V(\boldsymbol{\theta}^\circ) \stackrel{\text{def}}{=} \sup_{\boldsymbol{\theta}, \boldsymbol{\theta}' \in \mathcal{B}'(\epsilon, \boldsymbol{\theta}^\circ)} \sup_{\gamma \in S^p} \frac{\gamma^\top V(\boldsymbol{\theta})\gamma}{\gamma^\top V(\boldsymbol{\theta}')\gamma}.$$

Notice that under the condition $\mathfrak{A}_\epsilon V(\cdot) \leq \nu_1$, the topology induced by the metric $\mathfrak{S}(\cdot, \cdot)$ is (locally) equivalent to the Euclidean topology and the set $\mathcal{B}(\epsilon, \boldsymbol{\theta}^\circ)$ can be well approximated by the ellipsoid $\mathcal{B}'(\epsilon, \boldsymbol{\theta}^\circ)$ and computing the local entropy $\mathbb{Q}(\epsilon, \cdot)$ can be reduced to the Euclidean case; see Lemma 5.4 for more detail.

Now we are ready to state an exponential bound for the contrast process in the smooth case.

**Theorem 2.8.** *Assume that* (EG) *and* (ED) *hold true with some* $\nu_0$ *and* $\overline{\lambda} > 0$. *Suppose that there is a constant* $\epsilon > 0$ *such that* $\epsilon\rho/(1-\rho) \leq \overline{\lambda}$ *and for a fixed* $\nu_1 \geq 1$ *and each* $\boldsymbol{\theta} \in \Theta$, *it holds*

$$\mathfrak{A}_\epsilon V(\boldsymbol{\theta}) \leq \nu_1. \tag{2.12}$$

*Let for some* $\rho, s < 1$ *the function* $\mathfrak{M}_\epsilon(\boldsymbol{\theta}^\circ, \boldsymbol{\theta}_0) = \inf_{\boldsymbol{\theta} \in \mathcal{B}(\epsilon, \boldsymbol{\theta}^\circ)} \mathfrak{M}(\boldsymbol{\theta}, \boldsymbol{\theta}_0)$ *fulfill*

$$\mathfrak{H}_\epsilon(\rho, s) \stackrel{\text{def}}{=} \log\left(\omega_p^{-1} \epsilon^{-p} \int_\Theta \sqrt{\det V(\boldsymbol{\theta})} \exp\big\{-\rho(1-s)\mathfrak{M}_\epsilon(\boldsymbol{\theta}, \boldsymbol{\theta}_0)\big\} d\boldsymbol{\theta}\right) < \infty,$$

*where* $\omega_p$ *is the Lebesgue measure of the unit ball in* $\mathbb{R}^p$. *Then it holds*

$$\mathfrak{Q}(\rho, s) \leq (1-\rho)\mathbb{Q}_p + \frac{2\nu_0^2 \epsilon^2 \rho^2}{1-\rho} + 2p\log(\nu_1) + \mathfrak{H}_\epsilon(\rho, s).$$

**Remark 2.1.** The conditions of this theorem are very mild. $(EG)$ only requires that $L(\boldsymbol{\theta}, \boldsymbol{\theta}_0)$ has exponential moments. $(ED)$ requires a similar condition for the centered and normalized gradient $\nabla L(\boldsymbol{\theta})$. The inequalities (2.12) are equivalent to uniform continuity of the function $V(\boldsymbol{\theta})$.

**Remark 2.2.** The presented exponential bound requires that the value $\mathfrak{H}_\epsilon(\rho, s)$ is finite. Fortunately it can be easily checked in typical situations. A typical example is given in Section 3 which deals with the i.i.d. case.

### 2.4.1. *A risk bound for $\widetilde{\boldsymbol{\theta}} - \boldsymbol{\theta}_0$*

Our main result controls the risk of the minimum contrast estimator in terms of the rate function $\mathfrak{M}(\boldsymbol{\theta}, \boldsymbol{\theta}_0)$. In the case of the smooth contrast, this result may be used to bound the classical estimation loss $\widetilde{\boldsymbol{\theta}} - \boldsymbol{\theta}_0$. The idea is to bound from the rate function $\mathfrak{M}(\boldsymbol{\theta}, \boldsymbol{\theta}_0)$ by a quadratic function in a vicinity of the point $\boldsymbol{\theta}_0$ and next to make use of the concentration property of $\widetilde{\boldsymbol{\theta}}$.

Note that for any $\mu$, it obviously holds $\mathfrak{M}(\mu, \boldsymbol{\theta}_0, \boldsymbol{\theta}_0) = 0$ and a simple algebra yields for the gradient of $\mathfrak{M}(\mu, \boldsymbol{\theta}_0, \boldsymbol{\theta}_0)$

$$\begin{aligned} \nabla \mathfrak{M}(\mu, \boldsymbol{\theta}, \boldsymbol{\theta}_0)\big|_{\boldsymbol{\theta}=\boldsymbol{\theta}_0} &= \frac{d}{d\boldsymbol{\theta}} \mathfrak{M}(\mu, \boldsymbol{\theta}, \boldsymbol{\theta}_0)\big|_{\boldsymbol{\theta}=\boldsymbol{\theta}_0} \\ &= -\mu I\!\!E \nabla L(\boldsymbol{\theta})\big|_{\boldsymbol{\theta}=\boldsymbol{\theta}_0} = -\mu \nabla I\!\!E L(\boldsymbol{\theta}_0) = 0. \end{aligned}$$

So, $\mathfrak{M}(\mu, \boldsymbol{\theta}_0, \boldsymbol{\theta}_0)$ can be majorated from below and from above in a vicinity of $\boldsymbol{\theta}_0$ by the Taylor expansion of the second order. The same behavior can be expected for the optimized rate function $\mathfrak{M}(\boldsymbol{\theta}_0, \boldsymbol{\theta}_0)$. This argument and the concentration property from Corollary 2.5 lead to the following result:

**Corollary 2.9.** *Suppose the conditions of Theorem 2.8 are satisfied and also for some $r > 0$, the function $\mathfrak{M}(\boldsymbol{\theta}, \boldsymbol{\theta}_0)$ fulfills*

$$\mathfrak{M}(\boldsymbol{\theta}, \boldsymbol{\theta}_0) \geq (\boldsymbol{\theta} - \boldsymbol{\theta}_0)^\top V_0 (\boldsymbol{\theta} - \boldsymbol{\theta}_0), \qquad \boldsymbol{\theta} \in \mathcal{A}(r, \boldsymbol{\theta}_0),$$

*for some positive matrix $V_0$. Then for any $\rho, s < 1$ and $\mathfrak{z} > 0$*

$$I\!\!P\big(\|\sqrt{V_0}(\widetilde{\boldsymbol{\theta}} - \boldsymbol{\theta}_0)\|^2 > \mathfrak{z}\big) \leq \mathfrak{Q}(\rho, s) \exp\{-\rho s \min\{\mathfrak{z}, r\}\}.$$

*Proof.* It is obvious that

$$\begin{aligned} \{\|\sqrt{V_0}(\widetilde{\boldsymbol{\theta}} - \boldsymbol{\theta}_0)\|^2 > \mathfrak{z}\} &\subseteq \{\|\sqrt{V_0}(\widetilde{\boldsymbol{\theta}} - \boldsymbol{\theta}_0)\|^2 > \mathfrak{z}, \widetilde{\boldsymbol{\theta}} \in \mathcal{A}(r, \boldsymbol{\theta}_0)\} \cup \{\widetilde{\boldsymbol{\theta}} \notin \mathcal{A}(r, \boldsymbol{\theta}_0)\} \\ &\subseteq \{\mathfrak{M}(\widetilde{\boldsymbol{\theta}}, \boldsymbol{\theta}_0) > \mathfrak{z}, \widetilde{\boldsymbol{\theta}} \in \mathcal{A}(r, \boldsymbol{\theta}_0)\} \cup \{\widetilde{\boldsymbol{\theta}} \notin \mathcal{A}(r, \boldsymbol{\theta}_0)\} \\ &= \{\widetilde{\boldsymbol{\theta}} \notin \mathcal{A}(r \wedge \mathfrak{z}, \boldsymbol{\theta}_0)\} \end{aligned}$$

and the result follows from Corollary 2.7. $\square$

In the case of i.i.d. observations, the function $\mathfrak{M}(\mu, \boldsymbol{\theta}, \boldsymbol{\theta}_0)$ and hence the matrix $V_0$ are proportional to the sample size $n$ and the result of Corollary 2.9 automatically yields the root-n consistency of $\widetilde{\boldsymbol{\theta}}$; see Section 3 for more details.

## 3. Quasi MLE for i.i.d. data

Let $\boldsymbol{Y} = (Y_1, \ldots, Y_n)$ be an i.i.d. sample from a distribution $P$. By $\boldsymbol{P}$ we denote the joint distribution of $\boldsymbol{Y}$. Let also $\mathcal{P} = (P_{\boldsymbol{\theta}}, \boldsymbol{\theta} \in \Theta \subset \mathbb{R}^p)$ be a parametric family. In contrast to the standard parametric hypothesis which assumes that $P \in \mathcal{P}$, in this section, we focus on the quality of estimation in the case when the underlying measure $P$ does not necessarily belong to the parametric family $\mathcal{P}$. We will see that in this case the maximum likelihood method estimates the point $\boldsymbol{\theta}_0$, which minimizes some special distance between $P$ and $P_{\boldsymbol{\theta}}$ over $\boldsymbol{\theta} \in \Theta$.

In the rest of this section, the family $\mathcal{P}$ and the underlying measure $P$ are assumed to be dominated by a measure $P_0$. We denote by $p(y, \boldsymbol{\theta})$ and $p(y)$ the corresponding densities: $p(y, \boldsymbol{\theta}) = dP_{\boldsymbol{\theta}}/dP_0(y)$, $p(y) = dP/dP_0(y)$. The maximum likelihood estimator $\widetilde{\boldsymbol{\theta}}$ of the underlying parameter $\boldsymbol{\theta}_0$ is computed as follows:

$$\widetilde{\boldsymbol{\theta}} = \operatorname*{argmax}_{\boldsymbol{\theta} \in \Theta} L(\boldsymbol{\theta}) = \operatorname*{argmax}_{\boldsymbol{\theta} \in \Theta} \sum_{i=1}^n \ell(Y_i, \boldsymbol{\theta}),$$

where $\ell(Y, \boldsymbol{\theta}) = \log p(Y, \boldsymbol{\theta})$. Denote $\ell(Y, \boldsymbol{\theta}, \boldsymbol{\theta}') = \ell(Y, \boldsymbol{\theta}) - \ell(Y, \boldsymbol{\theta}')$ and

$$\mathfrak{m}(\mu, \boldsymbol{\theta}, \boldsymbol{\theta}_0) = -\log E \exp\{\mu \ell(Y, \boldsymbol{\theta}, \boldsymbol{\theta}_0)\},$$

The i.i.d. structure of the observations $\boldsymbol{Y}$ implies that

$$\mathfrak{M}(\mu, \boldsymbol{\theta}, \boldsymbol{\theta}_0) = n\, \mathfrak{m}(\mu, \boldsymbol{\theta}, \boldsymbol{\theta}_0).$$

This enables us to redefine the function $\mu^*(\boldsymbol{\theta})$ in terms of the function $\mathfrak{m}(\cdot, \boldsymbol{\theta}, \boldsymbol{\theta}_0)$ corresponding to the marginal distribution $P$:

$$\mu^*(\boldsymbol{\theta}) = \operatorname*{argmax}_{\mu} \mathfrak{m}(\mu, \boldsymbol{\theta}, \boldsymbol{\theta}_0)$$

and $\mu(\boldsymbol{\theta})$ can be interpreted as an approximation of $\mu^*(\boldsymbol{\theta})$. Denote also

$$\mathfrak{m}(\boldsymbol{\theta}, \boldsymbol{\theta}_0) = \mathfrak{m}(\mu(\boldsymbol{\theta}), \boldsymbol{\theta}, \boldsymbol{\theta}_0),$$

and for $\zeta_1(\boldsymbol{\theta}) = \mu(\boldsymbol{\theta})\{\ell(Y_1, \boldsymbol{\theta}, \boldsymbol{\theta}_0) - E\ell(Y_1, \boldsymbol{\theta}, \boldsymbol{\theta}_0)\}$ define

$$\begin{aligned} v(\boldsymbol{\theta}) &= E\, \nabla \zeta_1(\boldsymbol{\theta})[\nabla \zeta_1(\boldsymbol{\theta})]^\top, \\ h(\delta, \gamma; \boldsymbol{\theta}) &= \log E \exp\left\{2\delta \frac{\gamma^\top \nabla \zeta_1(\boldsymbol{\theta})}{\sqrt{\gamma^\top v(\boldsymbol{\theta})\gamma}}\right\}. \end{aligned}$$

Notice that if $P$ coincides with $P_{\boldsymbol{\theta}_0}$ and $\mu(\boldsymbol{\theta})$ is constant in a vicinity of $\boldsymbol{\theta}_0$, then $v(\boldsymbol{\theta}_0)$ is the standard Fisher information matrix. One can easily check that

$$h(0, \gamma; \boldsymbol{\theta}) = 0, \quad \left.\frac{\partial h(\delta, \gamma; \boldsymbol{\theta})}{\partial \delta}\right|_{\delta=0} = 0, \quad \left.\frac{\partial^2 h(\delta, \gamma; \boldsymbol{\theta})}{\partial^2 \delta}\right|_{\delta=0} = 4.$$

It follows from Lemma 5.8 that for any $\nu_0 > 1$ and $\boldsymbol{\theta} \in \Theta$ there exists $\overline{\delta}(\boldsymbol{\theta}, \nu_0) > 0$ such that $h(\delta, \gamma; \boldsymbol{\theta}) \leq 2\nu_0^2 \delta^2$ for all $\gamma \in S^p$ and $\delta \leq \overline{\delta}(\boldsymbol{\theta}, \nu_0)$. We assume a slightly stronger condition that $\overline{\delta}(\boldsymbol{\theta})$ can be taken the same for all $\boldsymbol{\theta}$, i.e.

$$\sup_{\boldsymbol{\theta} \in \Theta} \sup_{\gamma \in S^p} h(\delta, \gamma; \boldsymbol{\theta}) \leq 2\nu_0^2 \delta^2, \qquad \delta \leq \overline{\delta}. \tag{3.1}$$

In some cases, the matrix $v(\boldsymbol{\theta})$ should be replaced by its regularization $\overline{v}(\boldsymbol{\theta})$ to ensure this property, see Section 4.2 for an example.

Independence of the $Y_i$'s implies that $V(\boldsymbol{\theta}) \stackrel{\text{def}}{=} \text{Cov}\{\nabla \zeta(\boldsymbol{\theta})\} = nv(\boldsymbol{\theta})$ and

$$H(\lambda, \gamma, \boldsymbol{\theta}) \stackrel{\text{def}}{=} \log I\!\!E \exp\Big\{2\lambda \frac{\gamma^\top \nabla \zeta(\boldsymbol{\theta})}{\sqrt{\gamma^\top V(\boldsymbol{\theta})\gamma}}\Big\} = nh(n^{-1/2}\lambda, \gamma; \boldsymbol{\theta})$$

for any $\lambda$ and any $\gamma \in S^p$. Therefore, if $n^{-1/2}\lambda \leq \overline{\delta}$, then by (3.1):

$$H(\lambda, \gamma, \boldsymbol{\theta}) \leq 2\nu_0^2 \lambda^2$$

and the condition $(ED)$ is fulfilled with $\overline{\lambda} \leq n^{1/2}\overline{\delta}$. Now one can easily reformulate Theorem 2.8 in terms of the marginal distribution $P$.

**Theorem 3.1.** *Assume (3.1) for some $\overline{\delta} > 0$ and $\nu_0 \geq 1$. Suppose that there are constants $\epsilon > 0$ and $\nu_1 \geq 1$ such that for each $\boldsymbol{\theta} \in \Theta$*

$$\mathfrak{A}_\epsilon v(\boldsymbol{\theta}) \leq \nu_1. \tag{3.2}$$

*Let also for some $s, \rho < 1$ such that $\epsilon \rho/(1-\rho) \leq n^{1/2}\overline{\delta}$*

$$\mathfrak{H}_\epsilon(\rho, s) \stackrel{\text{def}}{=} \log \left( \omega_p^{-1} \epsilon^{-p} \int_\Theta \sqrt{\det\{nv(\boldsymbol{\theta})\}} \exp\{-\rho(1-s)n\,\mathfrak{m}_\epsilon(\boldsymbol{\theta}, \boldsymbol{\theta}_0)\} d\boldsymbol{\theta} \right) < \infty,$$

*where $\mathfrak{m}_\epsilon(\boldsymbol{\theta}, \boldsymbol{\theta}_0) = \inf_{\boldsymbol{\theta}' \in \mathcal{B}(\epsilon, \boldsymbol{\theta})} \mathfrak{m}(\boldsymbol{\theta}, \boldsymbol{\theta}_0)$. Then the value $\mathfrak{Q}(\rho, s)$ from (2.2) fulfills*

$$\log \mathfrak{Q}(\rho, s) \leq (1-\rho)\mathbb{Q}_p + \frac{2\nu_0^2 \epsilon^2 \rho^2}{1-\rho} + 2p\log(\nu_1) + \mathfrak{H}_\epsilon(\rho, s).$$

The integral in $\mathfrak{H}_\epsilon(\rho, s)$ can be easily bounded in typical situations. The result presented below involves some conditions on the marginal rate function $\mathfrak{m}(\boldsymbol{\theta}, \boldsymbol{\theta}_0)$. Namely, it is assumed that this function is bounded from below by a quadratic polynom in a vicinity $\mathcal{A}_1(r, \boldsymbol{\theta}_0) \stackrel{\text{def}}{=} \{\boldsymbol{\theta} : \mathfrak{m}(\boldsymbol{\theta}, \boldsymbol{\theta}_0) \leq r\}$ of the point $\boldsymbol{\theta}_0$ for some fixed $r > 0$ and it increases at least logarithmically with the norm $\|\boldsymbol{\theta} - \boldsymbol{\theta}_0\|$ outside of this neighborhood.

In particularly, it is shown in Section 5 that for $n$ sufficiently large

$$\mathfrak{H}_\epsilon(\rho, s) \approx \log\left(1 + \frac{\omega_p^{-1} \pi^p}{|\mathfrak{a}_r^2 \epsilon^2 \rho(1-s)|^{p/2}}\right). \tag{3.3}$$

**Theorem 3.2.** *Assume (3.1) and let $\rho$ fulfill $\rho/(1-\rho) \leq n\overline{\delta}^2$. Suppose that (3.2) holds with $\epsilon = \sqrt{(1-\rho)/\rho}$. Let for some $r > 0$, there are a positive matrix $v_0$ and a constant $\mathfrak{a}_r > 0$ such that*

$$v(\boldsymbol{\theta}) \leq v_0, \qquad \mathfrak{m}(\boldsymbol{\theta},\boldsymbol{\theta}_0) \geq \mathfrak{a}_r^2(\boldsymbol{\theta}-\boldsymbol{\theta}_0)^\top v_0(\boldsymbol{\theta}-\boldsymbol{\theta}_0), \quad \forall \boldsymbol{\theta} \in \mathcal{A}_1(r,\boldsymbol{\theta}_0).$$

*Let for some $\beta > 0$, hold:*

$$C_r(\beta) \stackrel{\text{def}}{=} \int_{\Theta \setminus \mathcal{A}_1(r,\boldsymbol{\theta}_0)} \sqrt{\det\{v(\boldsymbol{\theta})\}} \exp\{-\beta \mathfrak{m}_\epsilon(\boldsymbol{\theta},\boldsymbol{\theta}_0)\} d\boldsymbol{\theta} < \infty.$$

*Finally, let $n$ be sufficiently large to ensure*

$$\mathfrak{b}_r(n) \stackrel{\text{def}}{=} \rho(1-s)nr - \beta r - \mathfrak{a}_r^{-1}\epsilon - (p/2)\log n \leq 0. \tag{3.4}$$

*Then for some $C$ depending on $\mathfrak{a}_r, \nu_0, \nu_1, C_r(\beta)$ only, it holds*

$$\log \mathfrak{Q}(\rho,s) \leq Cp + \frac{p}{2}\log\bigl(|(1-\rho)(1-s)|^{-1}\bigr),$$

This bound together with Corollary 2.9 yields

$$I\!P\bigl(n\mathfrak{a}_r^2 \|v_0^{1/2}(\widetilde{\boldsymbol{\theta}}-\boldsymbol{\theta}_0)\|^2 > \mathfrak{z} + pC(\rho,s)\bigr) \leq \exp\{-\rho s \min\{\mathfrak{z}, r\sqrt{n}\}\}$$

with $C(\rho,s) = C + \log\bigl(|(1-\rho)(1-s)|^{-1}\bigr)/2$. This result means root-n consistency of $\widetilde{\boldsymbol{\theta}}$ in a rather strong sense.

## 4. Examples

This section illustrates how the exponential bounds can be applied to some particular situations. To simplify technical details, we do not try to cover the most general case. Rather we aim to show that our basic conditions can be easily verified in typical situations.

### 4.1. Estimation in the exponential model

The exponential model assumes that the observations $\boldsymbol{Y} = (Y_1,\ldots,Y_n)$ are i.i.d. exponential random variables from the exponential law $P_\theta$ with an unknown parameter $\theta \in I\!R^+$: $P_\theta(Y_i > y) = \exp(-\theta y)$. In this example we focus on the classical parametric set-up assuming that the underlying measure $I\!P$ coincides with the product of $I\!P_{\theta_0}$ for some $\theta_0 \in I\!R^+$. The corresponding maximum likelihood contrast is given by

$$L(\theta) = \sum_{i=1}^n \ell(Y_i, \theta) = -\theta \sum_{i=1}^n Y_i + n\log(\theta)$$

yielding

$$\widetilde{\theta} = n \Big/ \sum_{i=1}^{n} Y_i, \qquad L(\widetilde{\theta}, \theta) = n\log(\widetilde{\theta}/\theta) + n(\theta/\widetilde{\theta} - 1) = n\mathcal{K}(\widetilde{\theta}, \theta),$$

where $\mathcal{K}(\theta, \theta') = \theta'/\theta - 1 - \log(\theta'/\theta)$ is the Kullback-Leibler divergence between the exponential laws $P_\theta$ and $P_{\theta'}$.

Define

$$h_1(\delta) \stackrel{\text{def}}{=} \log E \exp\{-\delta(\theta_0 Y_1 - 1)\}.$$

Then, with $u = \theta/\theta_0 - 1$, it holds

$$\mathfrak{m}(\mu, \theta, \theta_0) \stackrel{\text{def}}{=} -\log E_{\theta_0} \exp\{\mu\ell(\theta, \theta_0)\} = \mu[u - \log(1+u)] - h_1(\mu u).$$

Therefore, with

$$\begin{aligned} \mathfrak{m}_1^*(u) &= \max_{\mu}\{\mu[u - \log(1+u)] - h_1(\mu u)\}, \\ \mu_1^*(u) &= \operatorname*{argmax}_{\mu}\{\mu[u - \log(1+u)] - h_1(\mu u)\}, \end{aligned}$$

the optimal choice of $\mu(\theta)$ is given by $\mu^*(\theta) = \mu_1^*(u)$ leading to $\mathfrak{m}^*(\theta, \theta_0) = \mathfrak{m}_1^*(u)$ for $u = \theta/\theta_0 - 1$. For applying Theorem 3.1, we need a lower bound for $\mathfrak{m}_1^*(u)$. Simple algebra yields for $Y_1 \sim Exp(\theta_0)$

$$h_1(\delta) = \delta - \log(1+\delta), \qquad \mathfrak{m}(\mu, \theta, \theta_0) = \log(1+\mu u) - \mu\log(1+u),$$

so that

$$\mu_1^*(u) = \operatorname*{argmax}_{\mu}\{\log(1+\mu u) - \mu\log(1+u)\} = \frac{u - \log(1+u)}{u\log(1+u)}.$$

To simplify the calculations, we proceed further with the suboptimal choice $\mu(\theta) \equiv \mu = 1/2$ instead of $\mu^*(\theta) = \mu_1^*(u)$ leading to $\mathfrak{m}(\theta, \theta_0) \stackrel{\text{def}}{=} \mathfrak{m}(\mu, \theta, \theta_0) = \mathfrak{m}(u)$ with

$$\mathfrak{m}(u) \stackrel{\text{def}}{=} \log(1 + u/2) - 0.5\log(1+u) = \frac{1}{2}\log\left(1 + \frac{u^2}{4(1+u)}\right)$$

for $u = \theta/\theta_0 - 1 > -1$. It is easy to see that $\mathfrak{m}(u) \geq c_1 u^2$ for $|u| \leq 1$, and $\mathfrak{m}(u) \geq c_2 \log(1+u)$ for $u \geq 1$ with some $c_1, c_2 > 0$.

Next

$$\begin{aligned} \zeta_1(\theta) &\stackrel{\text{def}}{=} \mu\{\ell(Y_1, \theta) - E\ell(Y_1, \theta)\} = -\mu\theta(Y_1 - 1/\theta_0), \\ \nabla\zeta_1(\theta) &= -\mu(Y_1 - 1/\theta_0) \end{aligned}$$

so that with $\sigma^2 = \operatorname{Var} Y_1 = 1/\theta_0^2$ it holds $v(\theta) \stackrel{\text{def}}{=} E[\nabla\zeta_1(\theta)]^2 \equiv \mu^2\sigma^2 = 1/(4\theta_0^2)$,

$$\log E \exp\{\delta\nabla\zeta_1(\theta)/\sqrt{v(\theta)}\} \equiv h_1(\delta),$$

and the condition (3.1) is obviously satisfied with some $\nu_0^2 < \infty$. Similarly, the conditions (5.4) through (3.4) can be easily verified and Theorem 3.2 applied with $s = 0$ yields

$$I\!\!E \exp\{\rho L(\widetilde{\theta}, \theta_0)/2\} \equiv I\!\!E \exp\{\rho n \, \mathcal{K}(\widetilde{\theta}, \theta_0)/2\} \leq \frac{C}{(1-\rho)^{1/2}}. \quad (4.1)$$

An important feature of this result is that it applies for the unbounded and non-compact parameter set $(0, +\infty)$. Another corollary of (4.1) is that the true parameter $\theta_0$ is covered with a high probability by the confidence set $\mathcal{E}(\mathfrak{z})$ of the form

$$\mathcal{E}(\mathfrak{z}) = \{\theta \in \Theta : \theta/\widetilde{\theta} - 1 - \log(\theta/\widetilde{\theta}) \leq \mathfrak{z}/n\}$$

provided that $\mathfrak{z}$ is sufficiently large.

### 4.2. LAD contrast and median estimation

Median or more generally quantile estimation is known to be more robust and stable against outliers and it is frequently used in econometric studies; see Koenker (2005), Koenker and Xiao (2006).

Suppose we are given a sample $\boldsymbol{Y} = (Y_1, \ldots, Y_n)$. In the problem of median estimation, these random variables are assumed i.i.d. and we are interested in estimating the median $\theta_0$ which is a root of the equation

$$P(Y_1 \leq \theta_0) = P(Y_1 \geq \theta_0).$$

Alternatively, the median minimizes the value $E|Y_1 - \theta|$ provided that the expectation of $|Y_1|$ is finite. This remark leads to the natural estimator $\widetilde{\theta}$ of the median as the minimizer of the contrast $-L(\theta) = \sum_{i=1}^{n} |Y_i - \theta|$:

$$\widetilde{\theta} = \operatorname*{argmax}_\theta L(\theta) = \operatorname*{argmin}_\theta \sum_{i=1}^{n} |Y_i - \theta|.$$

If the $Y_i$'s are i.i.d. with the Laplace density $\exp(-|y - \theta_0|)/2$, then $L(\theta)$ coincides (up to a constant factor) with the log-likelihood. In the general case, $L(\theta)$ can be treated as a quasi log-likelihood contrast. Later we also briefly comment on the case when the $Y_i$'s are not i.i.d.

Assume first that $Y_i$ has the density $p_\theta(y) = p(y-\theta)$ where $p(\cdot)$ is a centrally symmetric function. To simplify the notation, we also assume that $\theta_0 = 0$. The general case can be reduced to this one by a simple change of variables. The density $p(y)$ is supposed to be positive and for $y > 0$ we define

$$\lambda(y) = -(2y)^{-1} \log[2P(Y_1 > y)].$$

Equivalently, we can write $P(Y > y) = \mathrm{e}^{-2y\lambda(y)}/2$ for $y \geq 0$. The case with $\lambda(y) \geq \lambda_0 > 0$ corresponds to light tails while $\lambda(y) \to 0$ as $|y| \to \infty$ means heavy tails of the distribution $P$. Below we focus on the most interesting case when $\lambda(y)$ is positive and monotonously decreases to zero in $y > 0$. For simplicity of presentation we also assume that $\lambda(y)$ is sufficiently regular and its first derivative $\lambda'(y)$ is uniformly continuous on $\mathbb{R}$. The assumption of heavy tails implies that $[y\lambda(y)]' \in [0,1]$ and hence,

$$|y\lambda'(y)| = |[y\lambda(y)]' - \lambda(y)| < 1.$$

Let

$$m(\theta) \stackrel{\text{def}}{=} E|Y_1 - \theta|, \qquad q(\theta) \stackrel{\text{def}}{=} P(Y_1 \leq \theta) - P(Y_1 > \theta).$$

Obviously $m'(\theta) \stackrel{\text{def}}{=} \partial m(\theta)/\partial \theta = q(\theta)$. It is also clear that $|q(\theta)| \leq 1$. Next, for $\theta \geq 0$, it holds

$$\ell'(y,\theta,\theta_0) \stackrel{\text{def}}{=} \frac{\partial}{\partial y}\ell(y,\theta,\theta_0) = \begin{cases} 0, & y \notin [0,\theta], \\ 2, & \text{otherwise,} \end{cases}$$

and $\ell(y,\theta,\theta_0) = -\theta$ for $y < 0$. Therefore, integration by parts yields

$$\begin{aligned}
Ee^{\mu\ell(Y_1,\theta,\theta_0)} &= -\int e^{\mu\ell(y,\theta,\theta_0)} dP(Y_1 > y) \\
&= e^{-\mu\theta} + \int \mu\ell'(y,\theta,\theta_0)e^{\mu\ell(y,\theta,\theta_0)}P(Y_1 > y)\, dy \\
&= e^{-\mu\theta} + 2\mu \int_0^\theta e^{\mu(2y-\theta)}P(Y_1 > y)\, dy \\
&= e^{-\mu\theta} + \mu e^{-\mu\theta}\int_0^\theta e^{2y[\mu - \lambda(y)]}\, dy
\end{aligned}$$

and similarly for $\theta < \theta_0$. We now fix $\mu(\theta) = \lambda(\theta)$. Monotonicity of $\lambda(y)$ implies

$$Ee^{\mu(\theta)\ell(Y_1,\theta,\theta_0)} = e^{-\theta\lambda(\theta)} + \lambda(\theta)e^{-\theta\lambda(\theta)}\int_0^\theta e^{2y[\lambda(\theta)-\lambda(y)]}dy \leq \{1+\theta\lambda(\theta)\}e^{-\theta\lambda(\theta)}.$$

Therefore, for $\theta > 0$,

$$\mathfrak{m}(\theta,\theta_0) \geq \theta\lambda(\theta) - \log\{1+\theta\lambda(\theta)\}. \tag{4.2}$$

The same low bound holds true for $\theta < 0$. For $\theta\lambda(\theta) \leq 1$ it obviously holds

$$\mathfrak{m}(\theta,\theta_0) \geq \theta^2\lambda^2(\theta)/2.$$

Now we check the condition (3.1). Define

$$\zeta_0(\theta) \stackrel{\text{def}}{=} E_0(|Y_1 - \theta| - |Y_1|) - (|Y_1 - \theta| - |Y_1|).$$

Then, for $\theta > 0$

$$\nabla \zeta_0(\theta) = \mathbf{1}(Y_1 \leq \theta) - \mathbf{1}(Y_1 > \theta) - q(\theta),$$
$$E|\nabla \zeta_0(\theta)|^2 = 1 - q^2(\theta),$$
$$\operatorname{Var} \zeta_0(\theta) = \operatorname{Var} \int_0^\theta \nabla \zeta_0(\theta) d\theta \leq \theta \int_0^\theta E|\nabla \zeta_0(\theta)|^2 d\theta = \theta \int_0^\theta \{1 - q^2(\theta)\} d\theta,$$

and

$$\theta^{-2} \operatorname{Var} \zeta_0(\theta) = \theta^{-1} \int_0^\theta \{1 - q^2(\theta)\} d\theta \to 0, \qquad \theta \to \infty$$

because $q(\theta) \to 1$. Next, $\zeta_1(\theta) \stackrel{\text{def}}{=} \lambda(\theta)\zeta_0(\theta)$ and

$$\nabla \zeta_1(\theta) = \partial \zeta_1(\theta)/\partial \theta = \lambda(\theta)\nabla \zeta_0(\theta) + \theta \lambda'(\theta)\zeta_0(\theta)/\theta.$$

Note that $|\nabla \zeta_0(\theta)| \leq 1$ and $|\zeta_0(\theta)/\theta| \leq 1$, and in addition $\lambda(\theta) \to 0$ and $\operatorname{Var}(\zeta_0(\theta)/\theta) \to 0$ as $\theta \to \infty$, while $|\theta \lambda'(\theta)|$ remains bounded by one. This easily implies the condition (3.1) for some fixed $\overline{\delta} > 0, \nu_0 \geq 1$, and $\overline{v}(\theta) \equiv 1$. Moreover, if $E|Y_1|^\gamma < \infty$ for some $\gamma > 0$, then the onditions of Theorem 3.2 are fulfilled. This theorem applied with $\rho = s$ and Corollary 2.4 lead to the bound for the loss $\widetilde{u} = |\widetilde{\theta} - \theta_0|$:

$$E \exp\{\rho^2 n [\widetilde{u}\lambda(\widetilde{u}) - \log\{1 + \widetilde{u}\lambda(\widetilde{u})\}]\} \leq \frac{C}{\rho^{1/2}(1-\rho)^{1/2}}$$

with some fixed constant $C$ provided that $n$ exceeds some minimal sample size $n_0$.

The case of independent but non i.i.d. observations can be again reduced to the considered case using $P = n^{-1} \sum_{i=1}^n P_i$ and defining the point $\theta_0$ as a root of the equation

$$\sum_{i=1}^n P_i(Y_i < \theta) = \sum_{i=1}^n P_i(Y_i > \theta).$$

### 4.3. Estimation of the location of a change point

Suppose the observations $\boldsymbol{Y} = (Y_1, \ldots, Y_n)$ follow the change point model:

$$Y_i = A \mathbf{1}(i \leq \theta) + \sigma \xi_i, \quad i = 1, \ldots, n, \tag{4.3}$$

where $\xi_i$ is a standard white Gaussian noise. Our goal is to estimate the change point location $\theta \in \Theta = \{1, \ldots, n-1\}$. The obtained results can be easily extended to the case of non-Gaussian errors under some exponential moment conditions.

We begin with the case when the amplitude $A$ is known. To estimate $\theta$, we use the maximum likelihood estimator

$$\widetilde{\theta}_A = \underset{\theta \in \Theta}{\operatorname{argmax}} L_A(\theta),$$

where the maximum likelihood contrast is given by

$$L_A(\theta) = \frac{A}{\sigma^2}\sum_{i=1}^{\theta} Y_i - \frac{A^2}{2\sigma^2}\theta = \frac{A^2}{\sigma^2}\min(\theta, \theta_0) - \frac{A^2\theta}{2\sigma^2} + \frac{A}{\sigma}\sum_{i=1}^{\theta}\xi_i.$$

Note that $L_A(\theta)$ is a Gaussian random variable for every $\theta$ with

$$M(\theta, \theta_0) \stackrel{\text{def}}{=} -\mathbb{E}L_A(\boldsymbol{\theta}) = \frac{A^2}{2\sigma^2}|\theta - \theta_0|,$$

$$D^2(\theta, \theta_0) \stackrel{\text{def}}{=} \operatorname{Var} L_A(\theta) = \frac{A^2}{\sigma^2}|\theta - \theta_0| = 2M(\theta, \theta_0).$$

This yields for any $\mu \geq 0$

$$\mathfrak{M}(\mu, \theta, \theta_0) = \mu M(\theta, \theta_0) - \mu^2 D^2(\theta, \theta_0)/2 = (\mu - \mu^2)M(\theta, \theta_0),$$

and the corresponding values $\mu^*(\theta), \mathfrak{M}^*(\theta, \theta_0)$ can be easily computed:

$$\mu^*(\theta) = 1/2, \qquad \mathfrak{M}^*(\theta, \theta_0) = M(\theta, \theta_0)/4.$$

Therefore, for $\rho < 1$, Theorem 2.1 implies

$$\mathbb{E}\exp\left\{\rho^2\frac{A^2}{4\sigma^2}|\widetilde{\theta} - \theta_0|\right\} \leq \sum_{\boldsymbol{\theta} \in \Theta}\exp\left\{-\frac{\rho(1-\rho)}{4}M(\theta, \theta_0)\right\}$$

$$\leq 2\sum_{k=0}^{\infty}\exp\left\{-\frac{\rho(1-\rho)A^2}{8\sigma^2}k\right\} = \frac{2}{1 - C(\rho)}$$

where $C(\rho) = \exp\{-\rho(1-\rho)A^2/(8\sigma^2)\}$. By Lemma 5.7

$$\mathbb{E}|\widetilde{\theta}_A - \theta_0|^r \leq C_1(r)\left(\sigma^2/A^2\right)^r$$

with some constant $C_1(r)$.

Now we switch to the case when $A > 0$ is an unknown parameter. In this case, we cannot use the contrast $L_A(\theta)$ because it strongly depends on $A$. To find a reasonable contrast, one can use the maximum likelihood principle. Considering $A$ as a nuisance parameter and maximizing $L_A(\theta)$ w.r.t. $A \geq 0$ leads to the following estimator:

$$\widetilde{\theta} = \underset{\theta}{\operatorname{argmax}}\left\{\max_{A \geq 0} L_A(\theta)\right\} = \underset{\theta}{\operatorname{argmax}}\frac{1}{2\sigma^2\theta}\left[\sum_{i=1}^{\theta} Y_i\right]_+^2,$$

where $[x]_+ = \max(x, 0)$. In what follows we deal with a slightly modified version of this estimator

$$\widetilde{\theta} = \operatorname*{argmax}_{\theta \in \Theta^n} L(\theta), \quad \text{with a new contrast} \quad L(\theta) = \frac{1}{\sigma\sqrt{\theta}} \sum_{i=1}^{\theta} Y_i,$$

which is again a Gaussian one. By the model equation (4.3), this contrast can be represented in the form:

$$L(\theta) = \frac{1}{\sqrt{\theta}} \sum_{i=1}^{\theta} \xi_i + \frac{A \min(\theta, \theta_0)}{\sigma\sqrt{\theta}}.$$

It is easy to see that the drift $M(\theta, \theta_0) = -\mathbb{E}L(\theta, \theta_0)$ satisfies

$$M(\theta, \theta_0) = \mathfrak{a}d(\theta, \theta_0)$$

with $\mathfrak{a} = \sigma^{-1}A\sqrt{\theta_0}$ and

$$d(\theta, \theta') = 1 - \sqrt{\min\{\theta/\theta', \theta'/\theta\}} = \begin{cases} 1 - \sqrt{\theta/\theta'}, & \theta \leq \theta', \\ 1 - \sqrt{\theta'/\theta}, & \theta \geq \theta'. \end{cases}$$

Similarly,

$$D^2(\theta, \theta') \stackrel{\text{def}}{=} \operatorname{Var} L(\theta, \theta') = \frac{2|\theta' - \theta|}{(\sqrt{\theta} + \sqrt{\theta'})\sqrt{\max(\theta, \theta')}} = 2d(\theta, \theta')$$

and obviously, $M(\theta, \theta_0) = \mathfrak{a}D^2(\theta, \theta_0)/2$. Also $D^2(\theta, \theta_0) \leq 2$ for all $\theta$. As $L(\theta)$ is a Gaussian contrast, it holds

$$\mu^*(\theta) = \frac{M(\theta, \theta_0)}{D^2(\theta, \theta_0)} = \frac{\mathfrak{a}}{2}, \qquad \mathfrak{M}^*(\theta, \theta_0) = \frac{\mathfrak{a}^2}{8}d(\theta, \theta_0);$$

see Example 1.1. Note that for every $\theta \in \Theta$, the value $\mathfrak{M}^*(\theta, \theta_0)$ is bounded by $\mathfrak{a}^2/8 = A^2\theta_0/(8\sigma^2)$. So, this example is quite special in the sense that the Kullback-Leibler divergence between measures $\mathbb{P}_{\theta_0}$ and $\mathbb{P}_\theta$ does not grow to infinity with $\theta$. We will see that this fact results in an extra loglog-factor in the bound for the minimum contrast.

For given $\epsilon > 0$ and $\theta^\circ \in \Theta$, the local ball $\mathcal{B}(\epsilon, \theta^\circ) = \{D(\theta, \theta^\circ) \leq \epsilon\}$ can be represented in the form

$$\mathcal{B}(\epsilon, \theta^\circ) = \{\theta : \theta^\circ(1 - \epsilon^2/2)^2 \leq \theta \leq \theta^\circ(1 - \epsilon^2/2)^{-2}\}.$$

and it can be transformed into the usual symmetric interval around $\log \theta^\circ$ by using the parameter $\log \theta$ instead of $\theta$:

$$\mathcal{B}(\epsilon, \theta^\circ) = \Big\{\theta : \big|\log \theta - \log \theta^\circ\big| \leq -2\log(1 - \epsilon^2/2)\Big\}.$$

This immediately implies that the local entropy $\mathbb{Q}(\epsilon, \theta^\circ)$ is bounded by $\overline{\mathbb{Q}} = 1$ for all $\theta^\circ \in \Theta$.

Let the measure $\pi(\cdot)$ assign the mass 1 to any point $\theta = 1, \ldots, n$. Then $\pi(\mathcal{B}(\epsilon, \theta^\circ))$ is equal to the number $\Pi_\epsilon(\theta)$ of points $\theta$ in $B(\epsilon, \theta^\circ)$, and it obviously holds $\Pi_\epsilon(\theta) \approx K(\epsilon)\theta$ with $K(\epsilon) = (1 - \epsilon^2/2)^{-2} - (1 - \epsilon^2/2)^2 \geq \epsilon^2$ for $\epsilon \leq 1$, so that (2.4) is fulfilled. Fix $\epsilon^2 = 1/2$. The trivial lower bound $\mathfrak{M}(\theta, \theta_0) \geq 0$ yields for $\mathfrak{H}_\epsilon(\rho, s)$ from (2.5) for any $s \leq 1$:

$$\mathfrak{H}_\epsilon(\rho, s) \leq \log\left(\sum_{\theta=1}^n \frac{1}{\Pi_\epsilon(\theta)}\right) \leq \log(C_1 \log n)$$

for some $C_1 > 0$. This yields by Theorem 2.3 and its Corollary 2.4 that

$$\mathbb{E} \exp\{\rho \mathfrak{a}^2 d(\widetilde{\theta}, \theta_0)/8\} \leq C_2 \log n. \tag{4.4}$$

Combining this with Lemma 5.7 yields

$$\mathbb{E}\left|\frac{A^2 \theta_0}{\sigma^2} d(\widetilde{\theta}, \theta_0)\right|^r \leq C |\log \log n|^r.$$

The extra $\log \log$-factor in this bound is due to the unbounded parameter set. In the case "classical" situation when the size $A$ of the jump is bounded away from zero and infinity and the true "relative" location $\theta_0/n$ is bounded away from the edge $0$ similar calculations (not presented here) lead to a bound $\mathbb{E} \exp\{C_1 \rho^2 A^2 |\widetilde{\theta} - \theta_0|\} \leq C_2$ which does not involve any extra log-term; see e.g. Csörgő and Horváth (1997) and references therein for asymptotic versions of this result.

It is also interesting to compare this result with the accuracy of the maximum likelihood method in the case, where the magnitude of jump $A$ is known. One can see that there is a payment for the adaptation to the nuisance parameter $A$ which is in form of an extra $\log \log$-factor. Another observation is that the accuracy of estimation strongly depends on the true location $\theta_0$, more precisely, on the value $\mathfrak{a}^2 = A^2 \theta_0 / \sigma^2$. In the "classical" situation this value is of order $n$ leading to the accuracy of order $n^{-1} \log \log(n)$. If the value $\mathfrak{a}^2$ is smaller in order than $n$, then the accuracy becomes worse by the same factor. In particular, if $A^2 \theta_0 / \sigma^2$ is of order one, then even consistency of $\widetilde{\theta}$ cannot be claimed.

## 5. Proofs

This section collects proofs of the main theorems and some auxiliary facts.

### 5.1. Proof of Theorem 2.2

Assume that $\boldsymbol{\theta}^\circ \in \Theta$. First we establish a local bound for the maximum of the process $L(\boldsymbol{\theta}, \boldsymbol{\theta}_0)$ over the local ball $\mathcal{B}(\epsilon, \boldsymbol{\theta}^\circ) = \{\boldsymbol{\theta} : \mathfrak{S}(\boldsymbol{\theta}, \boldsymbol{\theta}^\circ) \leq \epsilon\}$.

*Proof.* The main step of the proof is a bound for the stochastic component $\zeta(\boldsymbol{\theta}, \boldsymbol{\theta}^\sharp)$ over the ball $\mathcal{B}(\epsilon, \boldsymbol{\theta}^\circ)$ for a fixed $\boldsymbol{\theta}^\sharp \in \mathcal{B}(\epsilon, \boldsymbol{\theta}^\circ)$.

**Lemma 5.1.** *Assume that $\zeta(\boldsymbol{\theta})$ is a separable process satisfying for any given $\boldsymbol{\theta}^\circ \in \Theta$ the condition $(EL)$. Then for any given $\boldsymbol{\theta}^\sharp \in \mathcal{B}(\epsilon, \boldsymbol{\theta}^\circ)$ and any $\lambda$ with $\lambda/\epsilon \leq \overline{\lambda}$*

$$\log I\!\!E \exp\left\{\frac{\lambda}{\epsilon} \sup_{\boldsymbol{\theta} \in \mathcal{B}(\epsilon, \boldsymbol{\theta}^\circ)} \zeta(\boldsymbol{\theta}, \boldsymbol{\theta}^\sharp)\right\} \leq \mathbb{Q}(\epsilon, \boldsymbol{\theta}^\circ) + 2\nu_0^2 \lambda^2.$$

*Proof.* The proof is based on the standard chaining argument (see e.g. van der Vaart and Wellner (1996)). Without loss of generality, we assume that $\mathbb{Q}(\epsilon, \boldsymbol{\theta}^\circ) < \infty$. Then for any integer $k \geq 0$, there exists a $2^{-k}\epsilon$-net $\mathcal{D}_k(\epsilon, \boldsymbol{\theta}^\circ)$ in the local ball $\mathcal{B}(\epsilon, \boldsymbol{\theta}^\circ)$ having the cardinality $\mathbb{N}(2^{-k}\epsilon, \epsilon, \boldsymbol{\theta}^\circ)$. Using the nets $\mathcal{D}_k(\epsilon, \boldsymbol{\theta}^\circ)$ with $k = 1, \ldots, K-1$, one can construct a chain connecting an arbitrary point $\boldsymbol{\theta}$ in $\mathcal{D}_K(\epsilon, \boldsymbol{\theta}^\circ)$ and $\boldsymbol{\theta}^\sharp$. It means that one can find points $\boldsymbol{\theta}_k \in \mathcal{D}_k(\epsilon, \boldsymbol{\theta}^\circ)$, $k = 1, \ldots, K-1$, such that $\mathfrak{S}(\boldsymbol{\theta}_k, \boldsymbol{\theta}_{k-1}) \leq 2^{-k+1}\epsilon$ for $k = 1, \ldots, K$. Here $\boldsymbol{\theta}_K$ means $\boldsymbol{\theta}$ and $\boldsymbol{\theta}_0$ means $\boldsymbol{\theta}^\sharp$. Notice that $\boldsymbol{\theta}_k$ can be constructed recurrently: $\boldsymbol{\theta}_k = \tau_k(\boldsymbol{\theta}_{k+1})$, $k = K-1, \ldots, 1$, where

$$\tau_k(\boldsymbol{\theta}) = \operatorname*{argmin}_{\boldsymbol{\theta}' \in \mathcal{D}_k(\epsilon, \boldsymbol{\theta}^\circ)} \mathfrak{S}(\boldsymbol{\theta}, \boldsymbol{\theta}').$$

It obviously holds for $\boldsymbol{\theta} \in \mathcal{D}_K(\epsilon, \boldsymbol{\theta}^\circ)$

$$\zeta(\boldsymbol{\theta}, \boldsymbol{\theta}^\sharp) = \sum_{k=1}^K \zeta(\boldsymbol{\theta}_k, \boldsymbol{\theta}_{k-1}).$$

For $\xi(\boldsymbol{\theta}_k, \boldsymbol{\theta}_{k-1}) = \zeta(\boldsymbol{\theta}_k, \boldsymbol{\theta}_{k-1})/\mathfrak{S}(\boldsymbol{\theta}_k, \boldsymbol{\theta}_{k-1})$ it holds that

$$\zeta(\boldsymbol{\theta}_k, \boldsymbol{\theta}_{k-1}) = \mathfrak{S}(\boldsymbol{\theta}_k, \boldsymbol{\theta}_{k-1})\xi(\boldsymbol{\theta}_k, \boldsymbol{\theta}_{k-1}) = 2\epsilon\, c_k\, \xi(\boldsymbol{\theta}_k, \boldsymbol{\theta}_{k-1})$$

with $c_k = c_k(\boldsymbol{\theta}, \boldsymbol{\theta}^\circ) = \mathfrak{S}(\boldsymbol{\theta}_k, \boldsymbol{\theta}_{k-1})/(2\epsilon) \leq 2^{-k}$, and

$$\begin{aligned}
\sup_{\boldsymbol{\theta} \in \mathcal{D}_K(\epsilon, \boldsymbol{\theta}^\circ)} \zeta(\boldsymbol{\theta}, \boldsymbol{\theta}^\sharp) &\leq \sum_{k=1}^K \sup_{\boldsymbol{\theta}' \in \mathcal{D}_k(\epsilon, \boldsymbol{\theta}^\circ)} \zeta(\boldsymbol{\theta}', \tau_{k-1}(\boldsymbol{\theta}')) \\
&= 2\epsilon \sum_{k=1}^K \sup_{\boldsymbol{\theta}' \in \mathcal{D}_k(\epsilon, \boldsymbol{\theta}^\circ)} c_k \xi(\boldsymbol{\theta}', \tau_{k-1}(\boldsymbol{\theta}')).
\end{aligned}$$

Since $c_k \leq 2^{-k}$, Lemma 5.6 below and condition $(EL)$ imply

$$\log I\!\!E \exp\left\{\frac{\lambda}{\epsilon} \sup_{\boldsymbol{\theta} \in \mathcal{D}_K(\epsilon,\boldsymbol{\theta}^\circ)} \zeta(\boldsymbol{\theta},\boldsymbol{\theta}^\sharp)\right\}$$

$$\leq \log I\!\!E \exp\left\{2\lambda \sum_{k=1}^{K} \sup_{\boldsymbol{\theta}' \in \mathcal{D}_k(\epsilon,\boldsymbol{\theta}^\circ)} c_k \xi(\boldsymbol{\theta}', \tau_{k-1}(\boldsymbol{\theta}'))\right\}$$

$$\leq \sum_{k=1}^{K} 2^{-k} \log\left[I\!\!E \exp\left\{\sup_{\boldsymbol{\theta}' \in \mathcal{D}_k(\epsilon,\boldsymbol{\theta}^\circ)} 2^k c_k\, 2\lambda \xi(\boldsymbol{\theta}', \tau_{k-1}(\boldsymbol{\theta}'))\right\}\right]$$

$$\leq \sum_{k=1}^{K} 2^{-k} \log\left[\sum_{\boldsymbol{\theta}' \in \mathcal{D}_k(\epsilon,\boldsymbol{\theta}^\circ)} I\!\!E \exp\left\{2^k c_k\, 2\lambda \xi(\boldsymbol{\theta}', \tau_{k-1}(\boldsymbol{\theta}'))\right\}\right]$$

$$\leq \sum_{k=1}^{K} 2^{-k}\left\{\log \mathbb{N}(2^{-k}\epsilon, \epsilon, \boldsymbol{\theta}^\circ) + 2\nu_0^2 \lambda^2\right\}.$$

These inequalities with the separability of $\zeta(\boldsymbol{\theta},\boldsymbol{\theta}^\sharp)$ yield

$$\log I\!\!E \exp\left\{\frac{\lambda}{\epsilon} \sup_{\boldsymbol{\theta} \in \mathcal{B}(\epsilon,\boldsymbol{\theta}^\circ)} \zeta(\boldsymbol{\theta},\boldsymbol{\theta}^\sharp)\right\} = \lim_{K \to \infty} \log I\!\!E \exp\left\{\frac{\lambda}{\epsilon} \sup_{\boldsymbol{\theta} \in \mathcal{D}_K(\epsilon,\boldsymbol{\theta}^\circ)} \zeta(\boldsymbol{\theta},\boldsymbol{\theta}^\sharp)\right\}$$

$$\leq \sum_{k=1}^{\infty} 2^{-k}\left\{2\nu_0^2 \lambda^2 + \log \mathbb{N}(2^{-k}\epsilon, \epsilon, \boldsymbol{\theta}^\circ)\right\} \leq 2\nu_0^2 \lambda^2 + \mathbb{Q}(\epsilon, \boldsymbol{\theta}^\circ)$$

which completes the proof of the lemma. □

Now we are prepared to complete the proof of the theorem. Denote

$$\boldsymbol{\theta}^\sharp = \operatorname*{argmax}_{\boldsymbol{\theta} \in \mathcal{B}(\epsilon,\boldsymbol{\theta}^\circ)} \left\{\mu(\boldsymbol{\theta}) I\!\!E L(\boldsymbol{\theta},\boldsymbol{\theta}_0) + \mathfrak{M}(\boldsymbol{\theta},\boldsymbol{\theta}_0)\right\}.$$

It is clear that

$$\sup_{\boldsymbol{\theta} \in \mathcal{B}(\epsilon,\boldsymbol{\theta}^\circ)} \left\{\mu(\boldsymbol{\theta}) L(\boldsymbol{\theta},\boldsymbol{\theta}_0) + \mathfrak{M}(\boldsymbol{\theta},\boldsymbol{\theta}_0)\right\}$$

$$\leq \mu(\boldsymbol{\theta}^\sharp) L(\boldsymbol{\theta}^\sharp,\boldsymbol{\theta}_0) + \mathfrak{M}(\boldsymbol{\theta}^\sharp,\boldsymbol{\theta}_0) + \sup_{\boldsymbol{\theta} \in \mathcal{B}(\epsilon,\boldsymbol{\theta}^\circ)} \zeta(\boldsymbol{\theta},\boldsymbol{\theta}^\sharp).$$

This yields by the Hölder inequality and Lemma 5.1 with $\lambda = \epsilon \rho/(1-\rho)$ that

$$\log I\!\!E \exp\left\{\sup_{\boldsymbol{\theta} \in \mathcal{B}(\epsilon,\boldsymbol{\theta}^\circ)} \rho\left[\mu(\boldsymbol{\theta}) L(\boldsymbol{\theta},\boldsymbol{\theta}_0) + \mathfrak{M}(\boldsymbol{\theta},\boldsymbol{\theta}_0)\right]\right\}$$

$$\leq \log I\!\!E \exp\left\{\rho\left[\mu(\boldsymbol{\theta}^\sharp) L(\boldsymbol{\theta}^\sharp,\boldsymbol{\theta}_0) + \mathfrak{M}(\boldsymbol{\theta}^\sharp,\boldsymbol{\theta}_0)\right] + \rho \sup_{\boldsymbol{\theta} \in \mathcal{B}(\epsilon,\boldsymbol{\theta}^\circ)} \zeta(\boldsymbol{\theta},\boldsymbol{\theta}^\sharp)\right\}$$

$$\leq \rho \log I\!\!E \exp\left\{\mu(\boldsymbol{\theta}^\sharp) L(\boldsymbol{\theta}^\sharp,\boldsymbol{\theta}_0) + \mathfrak{M}(\boldsymbol{\theta}^\sharp,\boldsymbol{\theta}_0)\right\}$$

$$+ (1-\rho) \log I\!\!E \exp\left\{\frac{\rho}{1-\rho} \sup_{\boldsymbol{\theta} \in \mathcal{B}(\epsilon,\boldsymbol{\theta}^\circ)} \zeta(\boldsymbol{\theta},\boldsymbol{\theta}^\sharp)\right\}$$

$$\leq (1-\rho)\mathbb{Q}(\epsilon,\boldsymbol{\theta}^\circ) + (1-\rho) 2\nu_0^2 \left|\frac{\epsilon\rho}{1-\rho}\right|^2$$

and the result follows. □

### 5.2. Proof of Theorem 2.3

Theorem 2.2 implies a local bound for the process $\mu(\boldsymbol{\theta})L(\boldsymbol{\theta},\boldsymbol{\theta}_0)+\mathfrak{M}(\boldsymbol{\theta},\boldsymbol{\theta}_0)$ over any ball $\mathcal{B}(\epsilon,\boldsymbol{\theta}^\circ)$. To derive a global bound we apply the following general fact:

**Lemma 5.2.** *Let $f(\boldsymbol{\theta})$ be a nonnegative function on $\Theta \subset \mathbb{R}^p$ and let for every point $\boldsymbol{\theta} \in \Theta$ a vicinity $U(\boldsymbol{\theta})$ be fixed such that $\boldsymbol{\theta}' \in U(\boldsymbol{\theta})$ implies $\boldsymbol{\theta} \in U(\boldsymbol{\theta}')$. Let also a measure $\pi(U(\boldsymbol{\theta}))$ of the set $U(\boldsymbol{\theta})$ fulfill for every $\boldsymbol{\theta}^\circ \in \Theta$*

$$\sup_{\boldsymbol{\theta} \in U(\boldsymbol{\theta}^\circ)} \frac{\pi(U(\boldsymbol{\theta}))}{\pi(U(\boldsymbol{\theta}^\circ))} \leq \nu. \tag{5.1}$$

*Then*

$$\sup_{\boldsymbol{\theta} \in \Theta} f(\boldsymbol{\theta}) \leq \nu \int_\Theta f^*(\boldsymbol{\theta}) \frac{1}{\pi(U(\boldsymbol{\theta}))} d\pi(\boldsymbol{\theta})$$

*with*

$$f^*(\boldsymbol{\theta}) \stackrel{\text{def}}{=} \sup_{\boldsymbol{\theta}' \in U(\boldsymbol{\theta})} f(\boldsymbol{\theta}').$$

*Proof.* For every $\boldsymbol{\theta}^\circ \in \Theta$

$$\int_\Theta f^*(\boldsymbol{\theta}) \frac{1}{\pi(U(\boldsymbol{\theta}))} d\pi(\boldsymbol{\theta}) \geq \int_{U(\boldsymbol{\theta}^\circ)} f^*(\boldsymbol{\theta}) \frac{1}{\pi(U(\boldsymbol{\theta}))} d\pi(\boldsymbol{\theta})$$

$$\geq f(\boldsymbol{\theta}^\circ) \int_{U(\boldsymbol{\theta}^\circ)} \frac{1}{\pi(U(\boldsymbol{\theta}))} d\pi(\boldsymbol{\theta})$$

because $\boldsymbol{\theta} \in U(\boldsymbol{\theta}^\circ)$ implies $\boldsymbol{\theta}^\circ \in U(\boldsymbol{\theta})$ and hence, $f(\boldsymbol{\theta}^\circ) \leq f^*(\boldsymbol{\theta})$. Now by (5.1)

$$\int_\Theta f^*(\boldsymbol{\theta}) \frac{1}{\pi(U(\boldsymbol{\theta}))} d\pi(\boldsymbol{\theta}) \geq \frac{f(\boldsymbol{\theta}^\circ)}{\nu} \int_{U(\boldsymbol{\theta}^\circ)} \frac{1}{\pi(U(\boldsymbol{\theta}^\circ))} d\pi(\boldsymbol{\theta}) = f(\boldsymbol{\theta}^\circ)/\nu$$

as required. □

We are going to apply Lemma 5.2 with

$$f(\boldsymbol{\theta}) = \exp\{\rho[\mu(\boldsymbol{\theta})L(\boldsymbol{\theta},\boldsymbol{\theta}_0) + s\mathfrak{M}(\boldsymbol{\theta},\boldsymbol{\theta}_0)]\}.$$

In view of the definition of $\mathfrak{M}_\epsilon(\boldsymbol{\theta}^\circ,\boldsymbol{\theta}_0) = \min_{\boldsymbol{\theta} \in \mathcal{B}(\epsilon,\boldsymbol{\theta}^\circ)} \mathfrak{M}(\boldsymbol{\theta},\boldsymbol{\theta}_0)$ it follows from the local bound of Theorem 5.1 that

$$\log \mathbb{E} \exp\Big\{ \sup_{\boldsymbol{\theta} \in \mathcal{B}(\epsilon,\boldsymbol{\theta}^\circ)} \rho\big[\mu(\boldsymbol{\theta})L(\boldsymbol{\theta},\boldsymbol{\theta}_0) + s\mathfrak{M}(\boldsymbol{\theta},\boldsymbol{\theta}_0)\big]\Big\}$$

$$\leq -\rho(1-s)\mathfrak{M}_\epsilon(\boldsymbol{\theta}^\circ,\boldsymbol{\theta}_0) + (1-\rho)\mathbb{Q}(\epsilon,\boldsymbol{\theta}^\circ) + \frac{2\nu_0^2\epsilon^2\rho^2}{1-\rho}.$$

and the theorem follows directly from Lemma 5.2.

### 5.3. Proof of Theorems 2.8

Below by $C_p$ we denote a generic constant (not necessarily the same) which only depends on the dimensionality $p$. First we show that the differentiability condition $(ED)$ implies the local moment condition $(EL)$.

**Lemma 5.3.** *Assume that $(ED)$ holds with some $\nu_0$ and $\overline{\lambda}$. Then for any $\boldsymbol{\theta}, \boldsymbol{\theta}' \in \Theta$ and any $\lambda$ with $|\lambda| \leq \overline{\lambda}$,*

$$\log I\!\!E \exp\left\{2\lambda \frac{\zeta(\boldsymbol{\theta}, \boldsymbol{\theta}')}{\mathfrak{S}(\boldsymbol{\theta}, \boldsymbol{\theta}')}\right\} \leq 2\nu_0^2 \lambda^2. \tag{5.2}$$

*Proof.* For $\boldsymbol{\theta}, \boldsymbol{\theta}' \in \Theta$, denote $\boldsymbol{u} = \boldsymbol{\theta}' - \boldsymbol{\theta}$. With these notations

$$L(\boldsymbol{\theta}, \boldsymbol{\theta}') = \boldsymbol{u}^\top \int_0^1 \nabla L(\boldsymbol{\theta} + t\boldsymbol{u}) dt.$$

Similar expressions hold for $I\!\!E L(\boldsymbol{\theta}, \boldsymbol{\theta}')$ and for $\zeta(\boldsymbol{\theta}, \boldsymbol{\theta}') = L(\boldsymbol{\theta}, \boldsymbol{\theta}') - I\!\!E L(\boldsymbol{\theta}, \boldsymbol{\theta}')$:

$$\zeta(\boldsymbol{\theta}, \boldsymbol{\theta}') = \boldsymbol{u}^\top \int_0^1 \nabla \zeta(\boldsymbol{\theta} + t\boldsymbol{u}) dt.$$

The definition of $\mathfrak{S}(\boldsymbol{\theta}, \boldsymbol{\theta}')$ implies for any $t \in [0, 1]$

$$c(t) \stackrel{\text{def}}{=} \frac{\sqrt{\boldsymbol{u}^\top V(\boldsymbol{\theta} + t\boldsymbol{u})\boldsymbol{u}}}{\mathfrak{S}(\boldsymbol{\theta}, \boldsymbol{\theta}')} \leq 1,$$

and therefore Lemma 5.6 and (2.10) with $\gamma = \boldsymbol{u}/\|\boldsymbol{u}\|$ yield

$$\begin{aligned}
\log I\!\!E \exp\left\{2\lambda \frac{\zeta(\boldsymbol{\theta}, \boldsymbol{\theta}')}{\mathfrak{S}(\boldsymbol{\theta}, \boldsymbol{\theta}')}\right\} &= \log I\!\!E \exp\left\{2\lambda \int_0^1 c(t) \frac{\gamma^\top \nabla \zeta(\boldsymbol{\theta} + t\boldsymbol{u})}{\sqrt{\gamma^\top V(\boldsymbol{\theta} + t\boldsymbol{u})\gamma}} dt\right\} \\
&\leq \int_0^1 c(t) \log I\!\!E \exp\left\{2\lambda \frac{\gamma^\top \nabla \zeta(\boldsymbol{\theta} + t\boldsymbol{u})}{\sqrt{\gamma^\top V(\boldsymbol{\theta} + t\boldsymbol{u})\gamma}}\right\} dt \\
&\leq 2\nu_0^2 \lambda^2
\end{aligned}$$

as required. □

Due to the next lemma, the smoothness of the contrast implies that the topology induced by the metric $\mathfrak{S}(\cdot, \cdot)$ is locally equivalent to the Euclidean topology and computing the local entropy $\mathbb{Q}(\epsilon, \cdot)$ can be reduced to the Euclidean case. Recall the notation

$$\mathcal{B}'(\epsilon, \boldsymbol{\theta}^\circ) = \left\{\boldsymbol{\theta} : (\boldsymbol{\theta} - \boldsymbol{\theta}^\circ)^\top V(\boldsymbol{\theta}^\circ) (\boldsymbol{\theta} - \boldsymbol{\theta}^\circ) \leq \epsilon^2\right\}.$$

The definition of $\mathcal{B}(\epsilon, \boldsymbol{\theta})$ implies that $\mathcal{B}(\epsilon, \boldsymbol{\theta}^\circ) \subseteq \mathcal{B}'(\epsilon, \boldsymbol{\theta}^\circ)$.

**Lemma 5.4.** *Assume* $(ED)$ *with some* $\overline{\lambda}$, *and let, for some fixed* $\nu_1 \geq 1$, $\epsilon > 0$

$$\mathfrak{A}_\epsilon V(\boldsymbol{\theta}) \leq \nu_1, \qquad \boldsymbol{\theta} \in \Theta. \tag{5.3}$$

*Then*

- $(EL)$ *is fulfilled for* $\lambda \leq \overline{\lambda}$, *i.e. (5.2) holds for all* $\lambda \leq \overline{\lambda}$.
- $\sup_{\boldsymbol{\theta} \in \Theta} \mathbb{Q}(\epsilon, \boldsymbol{\theta}) \leq \mathbb{Q}_p + p \log(\nu_1)$, *where* $\mathbb{Q}_p$ *is the entropy of the unit ball in* $I\!\!R^p$ *in the Euclidean topology.*

*Proof.* The first claim is an immediate corollary of Lemma 5.3. Fix any $\boldsymbol{\theta}^\circ \in \Theta$. Linear transformation with the matrix $V^{-1}(\boldsymbol{\theta}^\circ)$ reduces the situation to the case when $V(\boldsymbol{\theta}^\circ) \equiv I$ and $\mathcal{B}'(\epsilon, \boldsymbol{\theta}^\circ)$ is a usual Euclidean ball for any $\epsilon_0 \leq \epsilon$. Moreover, by (5.3), each elliptic set $\mathcal{B}'(\epsilon_0, \boldsymbol{\theta})$ for $\boldsymbol{\theta} \in \mathcal{B}(\epsilon, \boldsymbol{\theta}^\circ)$ is nearly an Euclidean ball in the sense that the ratio of its largest and smallest axes (which is the ratio of the largest and smallest eigenvalues of $V^{-1}(\boldsymbol{\theta}^\circ) V^2(\boldsymbol{\theta}) V^{-1}(\boldsymbol{\theta}^\circ)$) is bounded by $\nu_1$. Therefore, for any $\epsilon_0 \leq \epsilon$, a Euclidean net $\mathcal{D}^e(\epsilon_0/\nu_1)$ with the step $\epsilon_0/\nu_1$ ensures a covering of $\mathcal{B}(\epsilon, \boldsymbol{\theta}^\circ)$ by the sets $\mathcal{B}(\epsilon_0, \boldsymbol{\theta}^\circ)$, $\boldsymbol{\theta}^\circ \in \mathcal{D}^e(\epsilon)$. Therefore, the corresponding covering number is bounded by $(\nu_1 \epsilon/\epsilon_0)^p$ yielding the claimed bound for the local entropy. □

Now we are ready to proceed with the proof of Theorem 2.8. We make use of the following technical result which helps to bound the global supremum of a random function over an integral of local maxima.

Consider the ellipsoid $\mathcal{B}'(\epsilon, \boldsymbol{\theta}^\circ) = \{\boldsymbol{\theta} : (\boldsymbol{\theta} - \boldsymbol{\theta}^\circ)^\top V(\boldsymbol{\theta}^\circ)(\boldsymbol{\theta} - \boldsymbol{\theta}^\circ) \leq \epsilon^2\}$. Its Lebesgue measure fulfills $\pi(\mathcal{B}'(\epsilon, \boldsymbol{\theta}^\circ)) = \omega_p \epsilon^p/\sqrt{\det\{V(\boldsymbol{\theta}^\circ)\}}$ where $\omega_p$ is the volume of the unit ball in $I\!\!R^p$. Condition (2.12) implies (5.1) with $\nu = \nu_1^p$ for $\pi(U(\boldsymbol{\theta})) = \pi(\mathcal{B}'(\epsilon, \boldsymbol{\theta}))$ and the Lebesgue measure $\pi$. Now the result follows from Theorem 2.3.

### 5.4. Proof of Theorem 3.2

We start with some technical lemmas.

**Lemma 5.5.** *Suppose that for some* $r > 0$, *there are a positive matrix* $v_0$ *and a constant* $\mathfrak{a}_r > 0$ *such that*

$$v(\boldsymbol{\theta}) \leq v_0, \qquad \mathfrak{m}(\boldsymbol{\theta}, \boldsymbol{\theta}_0) \geq \mathfrak{a}_r^2 (\boldsymbol{\theta} - \boldsymbol{\theta}_0)^\top v_0 (\boldsymbol{\theta} - \boldsymbol{\theta}_0), \quad \boldsymbol{\theta} \in \mathcal{A}_1(r, \boldsymbol{\theta}_0) \tag{5.4}$$

*Then for any* $\eta > 0$

$$\int_{\mathcal{A}_1(r, \boldsymbol{\theta}_0)} \sqrt{\det\{n v(\boldsymbol{\theta})\}} \exp\{-\eta\, n\, \mathfrak{m}_\epsilon(\boldsymbol{\theta}, \boldsymbol{\theta}_0)\} d\boldsymbol{\theta} \leq \mathfrak{a}_r^{-p} \big(\omega_p \epsilon^p + |\pi/\eta|^{p/2}\big).$$

*Proof.* The conditions of the lemma imply that for $\boldsymbol{\theta} \in \mathcal{A}_1(r, \boldsymbol{\theta}_0)$

$$\sqrt{n \mathfrak{m}_\epsilon(\boldsymbol{\theta}, \boldsymbol{\theta}_0)} \geq \big[\sqrt{n} \mathfrak{a}_r \|v_0^{1/2}(\boldsymbol{\theta} - \boldsymbol{\theta}_0)\| - \epsilon\big]_+.$$

Changing the variable $\boldsymbol{\theta}$ by $\boldsymbol{u} = (n\mathfrak{a}_r^2)^{1/2} v_0^{1/2}(\boldsymbol{\theta} - \boldsymbol{\theta}_0)$, yields in view of (5.4) that

$$\int_{\mathcal{A}_1(r,\boldsymbol{\theta}_0)} \exp\{-\eta\, n\mathfrak{m}_\epsilon(\boldsymbol{\theta},\boldsymbol{\theta}_0)\} \sqrt{\det\{nv(\boldsymbol{\theta})\}}\, d\boldsymbol{\theta}$$

$$\leq \frac{1}{\mathfrak{a}_r^p} \left( \int_{\|\boldsymbol{u}\|\leq\epsilon} d\boldsymbol{u} + \int_{\mathbb{R}^p} \exp\{-\eta\|\boldsymbol{u}\|^2\} d\boldsymbol{u} \right) \leq \mathfrak{a}_r^{-p}\bigl(\omega_p \epsilon^p + |\pi/\eta|^{p/2}\bigr)$$

as required. □

Next we bound the part of the integral $\mathfrak{H}_\epsilon(\rho, s)$ over the complement of $\mathcal{A}_1(r, \boldsymbol{\theta}_0)$. Namely, we aim to show that

$$\int_{\Theta \setminus \mathcal{A}_1(r,\boldsymbol{\theta}_0)} \sqrt{\det\{nv(\boldsymbol{\theta})\}} \exp\bigl\{-\rho(1-s)n\, \mathfrak{m}_\epsilon(\boldsymbol{\theta},\boldsymbol{\theta}_0)\bigr\} d\boldsymbol{\theta} \leq C_r(\beta) e^{-\mathfrak{b}_r(n)}. \quad (5.5)$$

Under (5.4), it obviuosly holds for $\boldsymbol{\theta} \in \Theta \setminus \mathcal{A}_1(r,\boldsymbol{\theta}_0)$ that $\mathfrak{m}_\epsilon(\boldsymbol{\theta},\boldsymbol{\theta}_0) \geq r - \mathfrak{a}_r^{-1}\epsilon/n$ and

$$\begin{aligned}
\rho(1-s)n\, \mathfrak{m}_\epsilon(\boldsymbol{\theta},\boldsymbol{\theta}_0) &\geq \beta \mathfrak{m}_\epsilon(\boldsymbol{\theta},\boldsymbol{\theta}_0) + \{\rho(1-s)n - \beta\}(r - \mathfrak{a}_r^{-1}\epsilon/n) \\
&\geq \beta \mathfrak{m}_\epsilon(\boldsymbol{\theta},\boldsymbol{\theta}_0) + \mathfrak{b}_r(n) + (p/2)\log n
\end{aligned}$$

and (5.5) follows by $\det\{nv(\boldsymbol{\theta})\} = n^p \det\{v(\boldsymbol{\theta})\}$.

Lemma 5.5 with $\eta = \rho(1-s)$, (5.5), and $\mathfrak{b}_r(n) \leq 0$ imply

$$\mathfrak{H}_\epsilon(\rho, s) \leq \mathfrak{a}_r^{-p}\left(1 + \frac{\omega_p^{-1}\pi^{p/2}}{|\epsilon^2 \rho(1-s)|^{p/2}}\right) + C_r(\beta)/(\omega_p \epsilon^p).$$

To finalize the proof, we apply Theorem 3.1 with $\epsilon$ defined by the equation $\epsilon^2 = (1-\rho)/\rho$.

$$\begin{aligned}
\log \mathfrak{Q}(\rho, s) &\leq (1-\rho)\mathbb{Q}_p + 2\nu_0^2 \rho + 2p\log(\nu_1) \\
&\quad + \log\left(1 + \frac{\omega_p^{-1}\pi^p \mathfrak{a}_r^{-p}}{|(1-\rho)(1-s)|^{p/2}} + \frac{\omega_p^{-1}C_r(\beta)\rho^{p/2}}{(1-s)^{p/2}}\right) \\
&\leq Cp + \frac{p}{2}\log\bigl(|(1-\rho)(1-s)|^{-1}\bigr)
\end{aligned}$$

where $C$ is a constant whose value depends on $\mathfrak{a}_r$, $\nu_0, \nu_1$, and $C_r(\beta)$. It is also used that $\mathbb{Q}_p \leq Cp$ and $\log \omega_p^{-1} \leq Cp$.

### 5.5. Auxiliary facts

**Lemma 5.6.** *For any r.v.'s $\xi_k$ and any nonnegative $\lambda_k$ such that $\Lambda = \sum_k \lambda_k \leq 1$*

$$\log \mathbb{E} \exp\left(\sum_k \lambda_k \xi_k\right) \leq \sum_k \lambda_k \log \mathbb{E} e^{\xi_k}. \quad (5.6)$$

*Proof.* Convexity of $e^x$ and concavity of $x^\Lambda$ imply

$$\mathbb{E}\exp\left\{\frac{\Lambda}{\Lambda}\sum_k \lambda_k(\xi_k - \log \mathbb{E}e^{\xi_k})\right\} \leq \mathbb{E}^\Lambda \exp\left\{\frac{1}{\Lambda}\sum_k \lambda_k(\xi_k - \log \mathbb{E}e^{\xi_k})\right\}$$

$$\leq \left\{\frac{1}{\Lambda}\sum_k \lambda_k \mathbb{E}\exp(\xi_k - \log \mathbb{E}e^{\xi_k})\right\}^\Lambda = 1.$$

**Lemma 5.7.** *Let $\xi$ be a nonnegative random variable and $\varphi(\lambda) = \log \mathbb{E}\exp(\lambda\xi)$. Then for any $r > 0$*

$$\left(\mathbb{E}\xi^r\right)^{1/r} \leq \inf_{\lambda:\,\varphi(\lambda)\geq r} \lambda^{-1}\varphi(\lambda). \tag{5.7}$$

*In particular, if $\varphi(\lambda) \leq a + \sigma^2\lambda^2$ for some $a, \sigma \geq 0$, then*

$$\left(\mathbb{E}\xi^r\right)^{1/r} \leq 2\sigma\sqrt{\max\{a, r/2\}}. \tag{5.8}$$

*Proof.* Consider the following function

$$f(x) = \begin{cases} \log^r(x) & \text{for } x \geq e^r, \\ xr^r/e^r & \text{for } x \leq e^r. \end{cases}$$

A simple algebra reveals that for $x > e^r$

$$\begin{aligned} f'(x) &= rx^{-1}\log^{r-1}(x), \\ f''(x) &= r(r-1)x^{-2}\log^{r-2}(x) - rx^{-2}\log^{r-1}(x) \\ &= rx^{-2}[r - 1 - \log(x)]\log^{r-2}(x) < 0. \end{aligned}$$

Since the function $f(x)$ is linear for $x \leq e^r$, it is concave for all $x \geq 0$. It is also easy to check that $[\log(x)]_+^r \leq f(x)$, because for $x \leq e^r$, the function $f(x)$ coincides with the tangent of $\log^r(x)$ at $x = e^r$. Therefore,

$$x^r = \lambda^{-r}\log^r\left(e^{\lambda x}\right) \leq \lambda^{-r}f(e^{\lambda x})$$

and the Jensen inequality implies for any $\lambda \geq 0$

$$\mathbb{E}\xi^r \leq \lambda^{-r}\mathbb{E}f(e^{\lambda\xi}) \leq \lambda^{-r}f(\mathbb{E}e^{\lambda\xi}) = \lambda^{-r}f\left(e^{\varphi(\lambda)}\right). \tag{5.9}$$

If $\varphi(\lambda) \geq r$, then $f(e^{\varphi(\lambda)}) = \log^r(e^{\varphi(\lambda)}) = \varphi^r(\lambda)$ and (5.7) follows from (5.9).

To prove (5.8), it remains to notice that the monotonicity of $f(\cdot)$ implies, in view of (5.9), that

$$\begin{aligned} \left(\mathbb{E}\xi^r\right)^{1/r} &\leq \inf_{\lambda:\,a+\sigma^2\lambda^2\geq r}\left\{\frac{a}{\lambda} + \sigma^2\lambda\right\} = \begin{cases} \sigma r(r-a)^{-1/2}, & a < r/2 \\ 2\sigma\sqrt{a}, & a \geq r/2 \end{cases} \\ &\leq \begin{cases} 2\sigma\sqrt{r/2}, & a < r/2 \\ 2\sigma\sqrt{a}, & a \geq r/2 \end{cases} \leq 2\sigma\sqrt{\max\{a, r/2\}}. \end{aligned}$$

**Lemma 5.8.** *Let a r.v. $\xi$ fulfill $\mathbb{E}\xi = 0$, $\mathbb{E}\xi^2 = 1$ and $\mathbb{E}\exp(\lambda_1|\xi|) = \varkappa < \infty$ for some $\lambda_1 > 0$. Then for any $\rho < 1$ there is a constant $C_1$ depending on $\varkappa$, $\lambda_1$ and $\rho$ only such that for $\lambda < \rho\lambda_1$*

$$\log \mathbb{E}e^{\lambda\xi} \leq C_1\lambda^2/2.$$

*Moreover, there is a constant $\lambda_2 > 0$ such that for all $\lambda \leq \lambda_2$*

$$\log \mathbb{E}e^{\lambda\xi} \geq \rho\lambda^2/2.$$

*Proof.* Define $h(x) = (\lambda - \lambda_1)x + m\log(x)$ for $m \geq 0$ and $\lambda < \lambda_1$. It is easy to see by a simple algebra that

$$\max_{x \geq 0} h(x) = -m + m\log\frac{m}{\lambda_1 - \lambda}.$$

Therefore for any $x \geq 0$

$$\lambda x + m\log(x) \leq \lambda_1 x + \log\left(\frac{m}{e(\lambda_1 - \lambda)}\right)^m.$$

This implies for all $\lambda < \lambda_1$

$$\mathbb{E}|\xi|^m \exp(\lambda|\xi|) \leq \left(\frac{m}{e(\lambda_1 - \lambda)}\right)^m \mathbb{E}\exp(\lambda_1|\xi|).$$

Suppose now that for some $\lambda_1 > 0$, it holds $\mathbb{E}\exp(\lambda_1|\xi|) = \varkappa(\lambda_1) < \infty$. Then the function $h_0(\lambda) = \mathbb{E}\exp(\lambda\xi)$ fulfills $h_0(0) = 1$, $h_0'(0) = \mathbb{E}\xi = 0$, $h_0''(0) = 1$ and for $\lambda < \lambda_1$,

$$h_0''(\lambda) = \mathbb{E}\xi^2 e^{\lambda\xi} \leq \mathbb{E}\xi^2 e^{\lambda|\xi|} \leq \frac{1}{(\lambda_1 - \lambda)^2}\mathbb{E}\exp(\lambda_1|\xi|).$$

This implies by the Taylor expansion for $\lambda < \rho\lambda_1$ that

$$h_0(\lambda) \leq 1 + C_1\lambda^2/2$$

with $C_1 = \varkappa(\lambda_1)/\{\lambda_1^2(1-\rho)^2\}$, and hence, $\log h_0(\lambda) \leq C_1\lambda^2/2$.